\documentclass[preprint,11pt,authoryear,3p]{article}

\usepackage{graphicx} 

\usepackage{verbatim}
\usepackage{multirow}
\usepackage{multicol}
\usepackage{lscape}
\usepackage{graphicx}
\usepackage{tocloft}
\usepackage{etoc}      
\usepackage{authblk}
\usepackage{algorithmicx}
\usepackage{algpseudocode}
\usepackage{booktabs}
\usepackage{colortbl}
\usepackage{subcaption}
\usepackage{algorithm}
\usepackage{amsmath}
\usepackage{amsfonts}
\usepackage{mathrsfs}
\usepackage[dvipsnames]{xcolor}
\usepackage{graphics}
\usepackage{tikz}
\usepackage{pgfplots}
\usepackage{url}
\usepgfplotslibrary{statistics}
\pgfplotsset{compat=newest}
\usepackage{titlesec}
\titleclass{\part}{top}
\titleformat{\part}[display]
  {\normalfont\huge\bfseries}{\partname\ \thepart}{20pt}{\Huge}
\titlespacing*{\part}{0pt}{50pt}{40pt}
\titleclass{\part}{straight}
\usetikzlibrary{graphs} 
\usetikzlibrary[graphs] 
\usetikzlibrary{quotes} 
\usetikzlibrary[quotes] 
\usetikzlibrary{arrows.meta} 
\usetikzlibrary{calc}
\definecolor{yellowV}{HTML}{F1C40F}
\definecolor{greenV}{HTML}{aceca1}
\definecolor{pinkV}{HTML}{8f2d56}
\definecolor{peachV}{HTML}{fec3a6}
\definecolor{bluV}{HTML}{73d2de}
\definecolor{peachDV}{HTML}{E4455F}
\definecolor{bluDV}{HTML}{227396}
\definecolor{greenDV}{HTML}{94C804}

\definecolor{lgreenlock}{HTML}{A8E57A}
\definecolor{violetlock}{HTML}{e0aaff}
\definecolor{dgreenlock}{HTML}{2a9d8f}
\definecolor{blulock}{HTML}{3a86ff}
\definecolor{pinklock}{HTML}{ff006e}
\definecolor{orangelock}{HTML}{FFA77B}

\allowdisplaybreaks

\title{Eco-Conscious Customers Behavior in Capacitated Two-Echelon Location-Routing Models for Sustainable Last-Mile Delivery}

\author[1]{Valentina Bonomi\footnote{Corresponding author.\\E-mail address: valentina.bonomi@tecnico.ulisboa.pt. Phone number: +39 331 463 2310.} }
\author[1]{Diana Jorge}
\author[1]{Tânia Ramos}
\author[1]{Ana Barbosa-Póvoa}
\affil[1]{CEGIST, Instituto Superior Técnico, Universidade de Lisboa, Av. Rovisco Pais, Lisboa,
            1049-001, 
            Portugal,
            Portugal}

\begin{document}

\maketitle

\begin{abstract}
This paper introduces a novel capacitated Two-Echelon Location–Routing Problem with Eco-conscious Customer Behavior (2E-LRP-ECB) aimed at enhancing the environmental sustainability of last-mile delivery (LMD) operations. The model jointly optimizes dynamic satellites location, vehicle routing, and customer delivery modes, explicitly accounting for (i) heterogeneous customer travel behaviors, (ii) heterogeneous fleet composition, and (iii) diverse emission profiles across both echelons. A piecewise linear formulation captures the additional emissions from first-echelon vehicle stops, while customer travel emissions are computed based on individual willingness and capacity to use zero-emission transport. The problem is solved exactly for a wide set of real-world-based instances under four operational strategies, differing in optimization objectives and second-echelon fleet composition. Computational experiments, including a case study with a major Portuguese LMD provider, highlight the environmental and operational trade-offs inherent to strategic and operational choices such as fleet composition, satellite activation, and customer pick-up policies. Results reveal that minimizing distance can lead to substantial increases in emissions, while emissions-oriented strategies leverage customer travel to achieve significant sustainability gains without compromising service efficiency. A multi-objective analysis using the $\epsilon$-constraint method produces Pareto frontiers and knee-point solutions, offering actionable insights for balancing operational efficiency and environmental impact in sustainable LMD design.
\end{abstract}

\section{Introduction} \label{sec:Intro}
The rapid expansion of e-commerce in recent years has fundamentally transformed the landscape of last-mile delivery (LMD). As the final stage of the delivery process, where packages travel from distribution centers to end customers, LMD faces substantial challenges directly impacting both competitiveness and sustainability. According to the European E-Commerce Report \cite{Lone2023}, Europe has experienced a 47\% increase in business-to-customer (B2C) e-commerce sales compared to 2019, accompanied by an 8\% rise in the number of internet users making online purchases within the same period. This accelerated shift toward online shopping has prompted new players to enter the urban delivery market, intensifying competition and elevating customer expectations \cite{Bonomi2025}. Current market demands now feature delivery time windows as narrow as two hours, enhanced reliability, and seamless return processes.
In response to these growing requirements, many retailers have turned to two-echelon distribution systems. Under this model, large vehicles transport goods from peripheral warehouses to intermediate satellites located near or within urban centers. Smaller, more agile vehicles then collect orders from these satellites to complete the last-mile delivery to end customers or orders are picked up from these satellites by customers. Although this approach can significantly reduce transit times and support timely fulfillment of customer orders, it also increases the number of commercial vehicles circulating in urban areas. This growth contributes to externalities such as elevated delivery costs, greater energy consumption, increased $CO_2$ emissions, noise pollution, and congestion. Indeed, within the logistics sector, LMD alone accounts for up to 30\% of total $CO_2$ emissions, highlighting the critical need for more sustainable transportation practices.
Despite extensive research on multi-echelon distribution and location-routing problems (\cite{Dellaert2019},\cite{BenMohamed2023},\cite{Lehmann2024}), relatively few studies have examined how customer behavior interacts with two-echelon structures, especially concerning environmental outcomes. While intermediate satellites are often praised for reducing vehicle travel within city centers, their potential benefits could be partially offset if customers choose to collect orders using private cars, thereby reintroducing traffic congestion and associated emissions.
In addition to fixed satellites, recent research has highlighted the role of \textit{dynamic satellites}, i.e., intermediate facilities that can be flexibly activated or relocated to better match daily demand and alleviate inner-city traffic (\cite{Sutrisno2023}, \cite{Liu2023MPLP}, \cite{Stokkink2025MicroHubs}). 
In most of these studies, mobility refers to the physical relocation or time–space deployment of satellites and lockers. 
Moreover, existing literature (\cite{Grabenschweiger2021},\cite{Tilk2021},\cite{Wang2022}) frequently simplifies customer travel behavior, assuming uniform choices across the population and neglecting crucial factors such as the distance to pick-up points, availability and convenience of public transport, and the customers' environmental consciousness. These oversimplifications can lead to inaccurate assessments of the true environmental impacts of last-mile delivery configurations, potentially undermining the effectiveness of proposed delivery strategies.

This paper addresses these critical gaps by formulating a mixed integer linear programming (MILP) model of a two-echelon location routing problem with eco-conscious customers' behavior (2E-LRP-ECB). The model integrates heterogeneous customer behavior, diverse vehicle types and dynamic satellites.
In contrast with previous studies, our approach interprets mobility not as the physical relocation of satellites, but as the dynamic \textit{activation} of facilities. 
Specifically, in each planning horizon a subset of candidate satellites is selected and activated according to demand conditions, while the remaining candidates remain inactive. 
We examine the problem by solving exactly several benchmark instances under varying conditions, specifically exploring four distinct scenarios based on the minimization of $CO_2$ emissions or distances and differing vehicle capacities and environmental friendliness. Additionally, a multi-objective analysis is conducted to assess the trade-offs between emissions reduction and delivery distances.

\subsection*{Contributions}
This paper provides multiple contributions:
\begin{itemize}
\item[a.] It introduces a novel scenario of the 2E-LRP-ECB explicitly incorporating heterogeneous customer behaviors, dynamic satellites activation and multiple vehicle types across both echelons. The developed 2E-LRP-ECB captures realistic assumptions about customer travel choices, varying vehicle emission profiles, and the interaction between satellite pick-up and home-delivery, significantly enhancing the practical relevance and accuracy of environmental impact assessments. Moreover, the daily selection of satellites sites from a larger candidate set significantly enhances the practical relevance and accuracy of the problem.
\item[b.] Through an exact solution approach applied to benchmark instances, detailed scenario-based experiments are performed. These experiments show the trade-offs inherent in vehicle type selection, emissions reduction strategies, and diverse customer behaviors, yielding valuable managerial insights.
\item[c.] A comprehensive multi-objective analysis is conducted to examine the interplay between emissions minimization and travel distances. This analysis provides decision-makers with critical insights into optimizing two-echelon distribution strategies that align both operational efficiency and sustainability objectives.
\end{itemize}
The remainder of this paper is organized as follows. Section \ref{sec:literature} reviews the literature on multi-echelon location-routing problems with an emphasis on sustainable urban last-mile delivery (LMD). Section \ref{sec:problemDesc} introduces the problem statement, while Section \ref{sec:model} details the proposed methodology. The experimental design is outlined in Section \ref{sec:solution}. Section \ref{sec:casestudy} presents the case study and the generation of benchmark instances, followed by the analysis of computational results in Section \ref{sec: settings}. Section \ref{sec:manager} discusses the main managerial insights derived from the findings. Finally, Section \ref{sec: conclusions} summarizes the key results and provides recommendations for future research on environmentally sustainable LMD solutions.

\section{Related Literature}\label{sec:literature}
Our research focuses on addressing a Two-Echelon Location Routing Problem (2E-LRP) within the Last Mile Delivery (LMD) sector with a particular emphasis on environmental impact.
The 2E-LRP originates as an extension of the classical Location Routing Problem (LRP), a well-established modeling framework in supply chain management that jointly addresses facility location decisions and vehicle routing to meet customer demands (\cite{Tordecilla2022,Alamatsaz2021}). While the LRP effectively captures core distribution trade-offs, the increasing complexity of urban logistics, driven by the rapid growth of e-commerce, sustainability concerns, and congestion issues, has led to the adoption of multi-echelon distribution structures.
The 2E-LRP integrates location and routing decisions across two hierarchical levels: the first-echelon, where goods are transported from warehouses to intermediate facilities (e.g., urban distribution satellites or satellites), and the second-echelon, where goods are delivered either through direct delivery from these facilities using smaller or more sustainable vehicles, or by allowing customers to collect their orders directly at the satellites. Early studies on multi-echelon logistics date back to \cite{Jacobsen1980}, who introduced a two-echelon model for newspaper distribution. 
Over the last decade, several extensions of the 2E-LRP have been proposed to incorporate real-world complexities. An exhaustive survey on 2E-LRPs scenarios can be found in \cite{Cuda2015}.
One of the key evolutions in last-mile logistics is the systematic use of heterogeneous fleets across echelons, often combining heavy-duty vehicles with lightweight or automated modes in dense urban areas. In the two-echelon setting, \cite{Yu2020Autonomous2E} proposed urban deliveries where conventional trucks in the first-echelon are combined with autonomous delivery devices in the second-echelon, formulating a MILP and designing tailored metaheuristics to scale to realistic instances (up to 200 customers). Complementarily, \cite{Liu2021VanRobotTRE} studied an e-grocery system that couples depot-to-hub vans with delivery robots in a multi-objective 2E-LRP, showing that integrating low-emission automated vehicles can reduce environmental burdens with limited operational penalties. Heterogeneity also arises on the human-resources side: \cite{Yu2021LRPTWCOOD} introduced a 2E-LRP with time windows and occasional drivers, highlighting cost and service trade-offs when blending company fleets with ad hoc crowd resources. 
\cite{Dellaert2019} studied the 2E-LRP with time windows where in the second-echelon goods have to be delivered to customers within specific time windows. 
\cite{BenMohamed2023} introduced the two-echelon stochastic multi-period capacitated LRP  where a stochastic multi-period characterization of the planning horizon is considered, shaping the evolution of the uncertain demand and costs.
Beyond fixed satellites, several works investigate dynamic or reconfigurable intermediate facilities to increase responsiveness and reduce inner-city travel.
\cite{Sutrisno2023} proposed a 2E-LRP problem with mobile satellites (2E-LRP-MS) to enhance LMD efficiency. It introduces a clustering-based simultaneous neighborhood search heuristic method for optimizing satellite location and routing solutions, reducing operational costs significantly.  \cite{Liu2023MPLP} study last-mile routing with mobile satellites, optimizing not only delivery tours but also time–space locker deployments. In these two works, satellites are treated as temporary or reconfigurable facilities whose location and deployment in time–space are optimized. 
Moreover, \cite{Lehmann2024} introduced a new scenario of the 2E-LRP that combines time windows, mixed pick-up and delivery demand and multiple trips in the second-echelon. An efficient Adaptive Large Neighborhood Search framework is presented to solve instances up to 100 customers.
\cite{Stokkink2025MicroHubs} addressed the optimal location of micro-hubs in multi-modal last-mile systems , providing structural insights on hub density and placement that are directly applicable to mobile or semi-mobile satellite planning.
In contrast to these approaches, where mobility typically refers to relocating or dynamically deploying lockers and hubs across time and space, our work adopts a different perspective. 
We consider a set of candidate satellite locations distributed across the city and allow the model to decide, for each planning horizon, which subset of these facilities is activated. A satellite becomes active when it is served by a first-echelon vehicle, but it may still operate without any second-echelon vehicle assignment if it is activated exclusively as a customer pick-up point. 

In recent years the focus has shifted to the Green LRP (G-LRP), a version of LRPs in which particular attention is put on the environmental impact of the system. Studies indicate that LMD contributes approximately 30\% of carbon dioxide ($CO_2$) emissions (\cite{Savadogo2021}) within the logistics sector. In response to these challenges, the European Union has set ambitious environmental targets for 2030, aiming to reduce greenhouse gas emissions by at least 55\% compared to 1990 levels. 
\cite{Alamatsaz2021} solved a multi objective capacitated G-LRP with the objective of finding the best location of facilities and simultaneously
design routes to satisfy customers’ stochastic demand with minimum total operating costs and total emitted carbon
dioxide. The authors proposed a combination of progressive hedging and genetic algorithms to solve the problem. 
\cite{Wang2022} designed a LRP with pick-up stations where deliveries are performed by green vehicles. Their purpose is to satisfy the total demand minimizing the activation costs of stations and the routing costs of vehicles. A branch-and-price algorithm is proposed to solve the problem. 
In a 2-echelon setting, 
\cite{Pitakaso2020} minimized the total fuel consumption of a 2E-LRP considering both distance and road conditions in both echelons. A new variable neighborhood strategy adaptive search algorithm is introduced to solve the problem.
\cite{Hajghani2023} proposed a multi-objective 2E-LRP aiming to minimize economic, environmental and social responsibility aspects. 
Another fundamental aspect of 2E-LRPs is the possibility for customers to travel to satellites to pick-up their orders. This approach can alleviate last-mile congestion, reduce delivery costs, and increase operational flexibility by shifting some of the distribution effort to the customer side. \cite{Tilk2021} developed a vehicle routing problem with delivery options in which each customer indicates a preferred delivery option, e.g. locker pickup, or retail shop. These alternate options may come with different time windows. The objective of the company is to serve all the customers minimizing routes' cost while complying with constraints on the  time windows associated with delivery locations, the vehicle capacity, and the minimum service level for each customer. Additionally, \cite{dosSantos2022} presented a 2E-LRP in which customers are divided into home-delivery and self pick-up. Moreover, customers visiting a locker may work as occasional couriers performing a delivery to another customer for a monetary compensation. Other studies do not set a priori the subset of customers that will receive their orders at home but leave it as a model decision. \cite{Grabenschweiger2021} incorporated customer satisfaction into a Vehicle Routing Problem (VRP) by introducing heterogeneous locker stations. Their objective is to minimize both travel costs and compensation costs offered to customers who choose to collect their orders at the lockers.
\cite{Wang2022} propose a LRP with pick-up stations, focusing on minimizing the total cost of activating stations and the routing costs of green vehicles. In their model, the delivery company decides whether customers are served through direct delivery or by having them retrieve their orders at a satellite.
Finally, \cite{Bonomi22} emphasize the role of customer behavior in reducing the environmental impact of last-mile delivery by presenting a LRP where the delivery method is based on each customer’s \textit{eco-conscious} behavior. Specifically, a customer opts to pick up the order if their own transport-related carbon emissions are lower than those associated with company delivery.

Building on the existing literature, our work introduces several novel features to the 2E-LRP framework. 
Starting from the contributions of \cite{Bonomi22}, our study moves from a single-echelon to a two-echelon framework while maintaining the decision structure for deliveries to end customers. In particular, each customer in our model is characterized by a travel profile defining both the maximum distance they can cover and the portion of this distance that can be traveled at zero emissions, a feature not yet explored in the two-echelon setting. 
Furthermore, our second-echelon employs a heterogeneous fleet where vehicles differ in both capacity and emissions levels. This stands in contrast to much of the existing literature, which often assumes homogeneous vehicles or omits explicit emission considerations.
Another novel aspect of our work lies in the way satellite activation costs are evaluated. Rather than imposing a nominal or symbolic cost, we incorporate the actual emissions stemming from vehicle stops (\cite{emissions}), thus achieving a more accurate representation of the total environmental impact. This approach contrasts with several prior studies, which typically treat satellite activation as a fixed overhead without reflecting its real effect on carbon emissions.

The main characteristics of the analyzed literature are summarized in Table \ref{tab:literature} and Table \ref{tab:literature2}. 
In Table \ref{tab:literature}, the column \textit{Multi-Echelon} reports whether the problem setting involves more than one echelon (a ‘$\checkmark$’ entry). 
Column \textit{Warehouses} specifies whether the first echelon includes one or multiple depots (typically the main warehouses in last-mile delivery). 
In the column \textit{Customers pickup}, a ‘$\checkmark$’ indicates that some customers may collect their orders directly from satellites instead of pure home delivery. 
The designation \textit{(Model)} means that the optimization model itself decides which customers are assigned to pickup versus home delivery, whereas \textit{(A priori)} means that the subset of customers using pickup is fixed in advance and treated as an input to the model. 
Columns \textit{Fleet (Hom./Het.)} describe the fleet composition in the first (E1) and second (E2) echelon. A fleet is considered heterogeneous if vehicles differ in capacity and/or emission profiles, while homogeneous fleets consist of identical vehicles. 
Column \textit{Env. Obj.} indicates whether environmental objectives, such as minimizing carbon emissions or fuel consumption, are part of the optimization.

Table \ref{tab:literature2} reports the characteristics of the instances used. 
Column \textit{Instance Type} distinguishes between datasets taken from the literature (\textit{Benchmark}), instances based on real data (\textit{Real}), and instances specifically designed for the study (\textit{Created}). 
Columns \textit{\#Customers} and \textit{\#Satellites} report the size of the tested instances. Intervals (e.g., [30,100]) indicate a range of values, while curly braces (e.g., \{25,50\}) specify discrete instance sizes. When two notations are combined with a slash, the study considers multiple types of instances (e.g., Real/Created). 
Finally, column \textit{Approach} indicates the solution strategy adopted: \textit{Exact} (solved to optimality), \textit{Heur.} (heuristics or metaheuristics), or a hybrid such as \textit{Exact/Heur}. 

\begin{table}[htbp!]
    \centering
    \scriptsize
        \caption{Characteristics of the reviewed literature in Location-Routing Problems}
    \resizebox{\textwidth}{!}{
  \begin{tabular}{lccccccc}
\toprule
\textbf{Authors} & \textbf{Multi-Echelon} & \textbf{Warehouses} & \textbf{Customers pickup} & \multicolumn{2}{c}{\textbf{Fleet (Hom./Het.)}} & \textbf{Env. Obj.}\\
\cmidrule(rr){5-6}
 &  &  &  & \textbf{E1} & \textbf{E2} &  &    \\
\midrule
\cite{Dellaert2019}   & $\checkmark$ & Multiple & - & Hom. & Hom. & -  \\
\cite{Yu2020Autonomous2E} & $\checkmark$ & Single & - & Hom. & Het. & $\checkmark$ \\
\cite{Alamatsaz2021} & - & Multiple & - & Hom. & - & $\checkmark$ \\
\cite{Grabenschweiger2021} & - & Single & $\checkmark$ (Model) & Hom & - & -  \\
\cite{Liu2021VanRobotTRE} & $\checkmark$ & Multiple & - & Hom. & Het. & $\checkmark$ \\
\cite{Tilk2021} & - & Single & $\checkmark$ (A priori) & Hom. & - & -  \\
\cite{Yu2021LRPTWCOOD} & $\checkmark$ & Multiple & $\checkmark$ (A priori) & Hom. & Het. & - \\
\cite{Bonomi22} & - & Single & $\checkmark$ (Model) & Hom & - & $\checkmark$  \\
\cite{dosSantos2022} &  $\checkmark$ & Single & $\checkmark$ (A priori) & Het. & Hom. & -  \\
\cite{Tordecilla2022} & -            & Multiple & - & Hom. & -           & -  \\
\cite{Wang2022}       & -            & Single   & $\checkmark$ (Model) & Hom. & - & $\checkmark$ \\
\cite{BenMohamed2023} & $\checkmark$ & Multiple & - & Hom. & Hom. & -  \\
\cite{Hajghani2023}   & $\checkmark$ & Multiple & - & Hom. & Hom. & $\checkmark$  \\
\cite{Liu2023MPLP} & $\checkmark$ & Single & $\checkmark$ (Model)  & Hom. & Hom. & $\checkmark$ \\
\cite{Sutrisno2023}   & $\checkmark$ & Multiple & - & Hom. & Hom. & - \\
\cite{Lehmann2024}    & $\checkmark$ & Multiple & - & Hom. & Hom. & -  \\
\cite{Stokkink2025MicroHubs} & $\checkmark$ & Multiple & - & Hom. & Hom. & $\checkmark$ \\
\midrule
\textbf{Current Work} & $\checkmark$ & Single   & $\checkmark$ (Model) 
                     & Het. & Het. (Capacity/Emissions) & $\checkmark$  \\
\bottomrule
\end{tabular}
}
    \label{tab:literature}
\end{table}

\begin{table}[htbp!]
    \centering
    \scriptsize
        \caption{Characteristics of the reviewed literature in Location-Routing Problems -- Instances Characteristics} 
    \resizebox{\textwidth}{!}{
  \begin{tabular}{lrrrr}
\toprule
\textbf{Authors} & \textbf{Instance Type} & \textbf{\#Customers} & \textbf{\#Satellites} & \textbf{Approach} \\
\midrule
\cite{Dellaert2019}   & Created & [15,100] & \{3,4,5\}  & Exact \\
\cite{Yu2020Autonomous2E} & Created & [50,200] & \{5,10\} & Heur. \\
\cite{Alamatsaz2021} & Real/Created & 36/[5,75] & \{3,4,5\} & Heur.\\
\cite{Grabenschweiger2021} & Created & [10,25] & [1,5]  &  Heur. \\
\cite{Liu2021VanRobotTRE} & Created/Real & [40,150] & [5,15] & Heur. \\
\cite{Tilk2021} & Created & \{25,50\} & \{5,10\}  & Exact\\
\cite{Yu2021LRPTWCOOD} & Created & [30,100] & \{3,5\} & Heur. \\
\cite{Bonomi22} & Created & \{50,100\} & \{5,10,15\} & Exact \\
\cite{dosSantos2022} & Benchmark & [45,100]& \{3,4,5\}& Exact \\
\cite{Tordecilla2022} & Created/Benchmark & \{8,10\}/[30,200] & \{2,3\}/[5,10]& Exact/Heur. \\
\cite{Wang2022} & Created & [21,100] & [3,30] & Exact \\
\cite{BenMohamed2023} & Created & \{15,20,50\}& 8 & Exact \\
\cite{Hajghani2023}   & Benchmark & [8,200] & [3,20]& Heur. \\
\cite{Liu2023MPLP} & Created & [50,200] & [5,10] & Heur. \\
\cite{Sutrisno2023}   & Created & [20,200]& \{5,10\}& Heur. \\
\cite{Lehmann2024}    & Created & 15/[30,100] & \{3,4,5\} & Exact/Heur. \\
\cite{Stokkink2025MicroHubs} & Real/Created & [10,300] & [2,10] & Heur. \\
\midrule
\textbf{Current Work} & Created & \{40,50,75,100,150\} & \{10\}& Exact \\
\bottomrule
\end{tabular}
}
    \label{tab:literature2}
\end{table}

\section{Problem description}\label{sec:problemDesc}
\begin{figure}[htbp!]
    \centering
    \scalebox{0.65}{
    \includegraphics[width=0.5\linewidth]{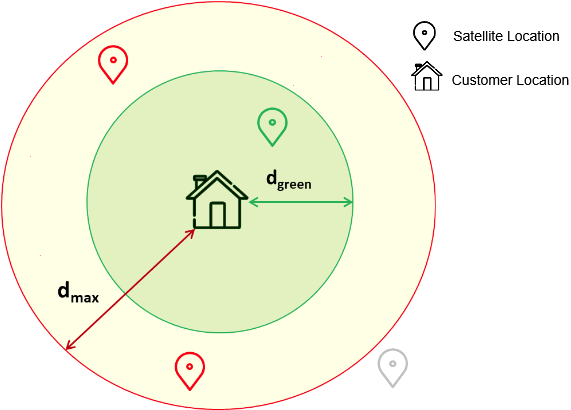}}
    \caption{Visual representation of $d_{max}$ and $d_{green}$ distances}
    \label{fig:dmax}
\end{figure}

We study a Two-Echelon Location–Routing Problem with Eco-conscious Customer Behavior (2E-LRP-ECB) that integrates location, assignment, and routing decisions with the aim of improving the environmental performance of last-mile distribution. The logistics system is structured in two echelons. In the first echelon, high-capacity combustion vehicles transport goods from a central warehouse to intermediate satellites acting as transshipment points. Unlike classical two-echelon settings, these satellites can be considered dynamic: the model selects the subset of candidate satellites to be activated according to demand distribution over the planning horizon. 
In the second echelon, activated satellites serve customers through two possible service modes: home delivery using small, possibly zero-emission vehicles (e.g., electric vans, bicycles, walking couriers), or customer pickup at a satellite. The assignment between these two modes is a model decision. Customers do not choose directly the delivery mode, the model determines whether each request is fulfilled by home delivery or by pickup, according to the optimization objectives.
The fleet is heterogeneous across echelons meaning that vehicles differ in capacity and emission profiles. In the second echelon, vehicles are not pre-assigned to satellites. Instead, the model decides the allocation of vehicle types and the construction of feasible routes departing from active satellites.  A satellite is considered active when it is supplied by at least a first-echelon vehicle, even if no second-echelon vehicle is dispatched from it, as customers may still be assigned to pick up their parcels directly at that location.
Each customer is characterized by two distance thresholds: (i) a maximum distance ($d_{max}$) within which pickup is admissible, and (ii) a shorter distance ($d_{green} \leq d_{max}$) defining whether the trip can be performed with a zero-emission mode such as walking or cycling. Package size further restricts admissibility, as very small parcels may be left in mailboxes without requiring pickup, whereas oversized items cannot be collected at satellites. These parameters jointly determine the conditions under which pickup can occur and the environmental impact of the chosen service mode.
Figure \ref{fig:dmax} illustrates how the distance thresholds $d_{max}$ and $d_{green}$ operate in the pickup setting. The inner green circle represents the range $d_{green}$, within which customers can access satellites using zero-emission modes such as walking or cycling. 
The larger red circle corresponds to $d_{max}$, defining the maximum admissible distance for pickup using emission-based transport modes. Satellites inside this range are eligible for pickup (colored in green and red) while the ones located outside $d_{max}$ are not admissible (colored in grey).

Environmental considerations are also embedded in the first-echelon supply of satellites. We explicitly account for the additional emissions arising from stops at satellites and along urban routes, where repeated acceleration and deceleration cycles increase fuel consumption. To capture this, emissions are modeled through a piecewise-linear function that scales with the number of stops and operating conditions.

\begin{figure}[htbp!]
\centering
\includegraphics[width=0.5\linewidth]{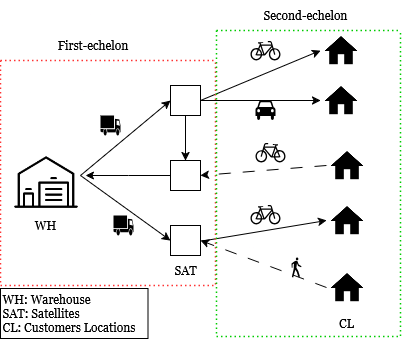}
\caption{Two-Echelon Location–Routing with dynamic satellites and mixed service modes (home delivery and pickup).}
\label{fig:2e}
\end{figure}

In conclusion, the main optimization decisions include: (i) activation of dynamic satellites, (ii) allocation of customers to satellites, (iii) selection of the service mode (home delivery or pickup) for each customer, (iv) routing of both supply and home-delivery vehicles, and (v) assignment of second-echelon vehicles to routes. Figure \ref{fig:2e} illustrates the considered setting.

\section{Methodology} \label{sec:model}
Formally, let $C = \{1,\dots,c_{max}\}$ be the set of customers that need a delivery and be $H = \{1,\dots,h_{max}\}$ the set of available intermediate satellites from which orders can be delivered or picked up by the customers. The set of vehicles used by the company is indicated as $K = K^1 \cup K^2$ where $K^1$ and $K^2$ contain first and second-echelon vehicles, respectively. Each vehicle is distinguished by a capacity $Q_k$, an emission factor $e_k$, a fuel consumption when loaded $\rho_k^*$ and a fuel consumption when empty $\rho_k$. Customer demand is indicated as $D_c$ and satellite capacity as $Q_h$.
The problem can be defined over a complete directed graph $G=(V,A)$ representing the road network where $V = \{0\} \cup C \cup H$ is the node set with node 0 representing the warehouse, and $A = A_1 \cup A_2 $ is the arc set with $A_1 = \{(i,j): i \in \{0\} \cup H, j \in H, i \neq j\}$ the arcs connecting the initial warehouse with the satellites and $A_2 = \{(i,j): i,j \in H \cup C, i \neq j\}$ the arcs connecting the satellites with the customers.
We indicate as $d_{ij}$ the non-negative distance between any two nodes $i, j \in V$ and we assume that distances satisfy the triangular inequality.  To separately account for arcs traveled by customers, we define as $H_c \subseteq H$ the subset of potential satellites located at a distance lower than or equal to $d_{max}^c$ from customer $c$. In addition, we define as $d_{green}^c$ the maximum distance threshold below which we assume any customer is willing to reach with zero-emission means of transport. Customers' emissions on arc $(c,h)$ with $c \in C$ and $h \in H_c$ is equal to $e_{ch}$ if $ d_{green}^c < d_{ch} \leq d_{max}^c$; 0 otherwise. 
Finally, the set $\Omega$ denotes the intervals used in the piecewise approximation of emissions due to stops. Each interval $\omega \in \Omega$ corresponds to a consecutive range of stops, and within that range a constant emission rate is applied. For example, interval $\omega=1$ may represent the case where a vehicle performs between 2 and 4 stops, and the associated emission rate is $e_1$. This discretization allows the model to capture the fact that emissions increase with the number of stops, while keeping the emission function linear within each interval.
Table \ref{tab:sets} and Table \ref{tab:parameters} summarize the sets and parameters introduced. 
\begin{table}[htbp!]
\centering
\scriptsize
\caption{Sets.}
\label{tab:sets}
\begin{tabular}{ll}
\toprule
\textbf{Set} & \textbf{Definition} \\
\midrule
$\{0\}$ & Warehouse \\
$H$ & Set of possible satellites \\
$C$ & Set of customers \\
$\mathcal{E}$ & Set of echelons \\
$K^1$ & Set of first-echelon vehicles \\
$K^2$ & Set of second-echelon vehicles \\
$K = K^1 \cup K^2$ & Set of all vehicles \\
$H_c \subseteq H$ & Set of satellites available for customer $c$ \\
$\Omega$ & Set of intervals of the piece-wise function related to the number of stops\\
\bottomrule
\end{tabular}
\end{table}
\begin{table}[htbp!]
\centering
\scriptsize
\caption{Parameters.}
\label{tab:parameters}
\begin{tabular}{ll}
\toprule
\textbf{Parameter} & \textbf{Definition} \\
\midrule
$Q_k$ & Capacity of vehicle $k \in K$ \\
$D_c$ & Demand of customer $c \in C$ \\
$Q_h$ & Capacity of satellite $h \in H$ \\
$e_k$ & Emissions of vehicle $k \in K$ \\
$e_{ch}$ & Emissions of customer $c$ to travel to satellite $h$ \\
$e_\omega$ & Emissions value for interval $\omega \in \Omega$\\
$d_{green}^c$ & Maximum no-emission distance traveled by customer $c$ for order pick-up\\
$d_{max}^c$ & Maximum distance traveled by customer $c$ for order pick-up\\
$d_{ij}$ & Distance on arc $(i,j) \in A$\\
\bottomrule
\end{tabular}
\end{table}

We formulate the problem described in Section \ref{sec:problemDesc} as a MILP.
Let us define, for each arc $(i,j) \in A$ and for each vehicle $k \in K$ a binary variable $x_{ij}^k$ taking value 1 if arc $(i,j)$ is traversed by vehicle $k$ and a positive continuous variable $z_{ij}^k \geq 0$ describing the load carried by vehicle $k$ on arc $(i,j)$. For each set of nodes $S \subset V$, let $\delta^+(S) = \{(i,j) \in A: i \in  S, j \notin S\}$ and 
$\delta^-(S) = \{(i,j) \in A: i \notin  S, j \in S\}$  be 
the set of arcs leaving and entering set $S$, respectively, with $\delta^+(i)= \delta^+(\{i\})$  and $\delta^-(i)= \delta^-(\{i\})$. 
The activation of satellites is regulated by the binary variable $l_h$ taking value 1 if the satellite $h \in H$ is active. The assignments of a second-echelon vehicle $k \in K^2$ to a specific satellite $h \in H$ and a specific customer $c \in C$ are regulated by the binary variable $q_k^h$ taking value 1 if vehicle $k$ is assigned to satellite $h$ and the binary variable $s_{c}^k$ taking value 1 if vehicle $k$ will visit customer $c$, respectively.
Moreover, for each customer $c$ and each satellite $h$ we define a binary variable $y_{ch}$ taking value 1 if customer $c$ receives its package at home from satellite $h$ and a binary variable $w_{ch}$ taking value 1 if instead the customer will pick up its package from satellite $h$.
Finally, for each vehicle $k \in K^1$, in each interval $\omega \in \Omega$, the continuous variable $f_{k\omega}$ is defined, counting the number of stops in the interval.
The decision variables are summarized in Table \ref{tab:variables}.
\begin{table}[htbp!]
\centering
\scriptsize
\caption{Decision variables.}
\begin{tabular}{ll}
\toprule
\textbf{Variable} & \textbf{Definition} \\
\midrule
$x_{ij}^{k}$ & 1 if arc $(i,j) \in A$ is traversed by vehicle $k \in K$; 0 otherwise \\
$y_{ch}$ & 1 if customer $c \in C$ is served from satellite $h \in H$; 0 otherwise\\
$w_{ch}$ & 1 if customer $c \in C$ picks up orders from satellite $h \in H$; 0 otherwise\\
$q_h^k$ & 1 if vehicle $k \in K^2$ is assigned to satellite $h \in H$; 0 otherwise\\
$s_{c}^k$ & 1 if customer $c \in C$ is assigned to vehicle $k \in K^2$; 0 otherwise\\
$l_{h}$ & 1 if satellite $h \in H$ is active; 0 otherwise\\ 
$f_{k\omega}$ & Number of stops of vehicle $k \in K^1$ in interval $\omega \in \Omega$\\
$z_{ij}^{k}$ & Load carried by vehicle $k \in K$ when traversing arc $(i,j)$ \\
\bottomrule
\end{tabular}
\label{tab:variables}
\end{table}

\newpage
The MILP is formulated as follows: 
\begin{alignat}{3}
& \min \sum_{k \in K^1}\sum_{\omega \in \Omega}e_{\omega}f_{k\omega} + \sum_{k \in K}\sum_{(i,j) \in A}e_kd_{ij}FCR_{ij}^k + \sum_{c \in C}\sum_{h \in H_c}e_{ch}w_{ch} \label{con:obj}\\
&\text{Where:}\notag\\
& FCR_{ij}^k = \rho_k x_{ij}^k + \frac{\rho_k^*-\rho_k}{Q_k}z_{ij}^k \quad \forall k \in K, (i,j) \in A \label{eq:FCR}\\
&\text{Subject to:}\notag\\
&\text{Echelon 1:}\notag\\
& \sum_{(0,j) \in \delta^+(0)}x_{0j}^{k} = \sum_{(i,0) \in \delta^-(0)}x_{i0}^{k}  \quad \forall k \in K^1 \label{con:warehouseflow} \\
&\sum_{(i,0) \in \delta^-(0)}x_{i0}^{k} \leq 1 \quad \forall k \in K^1 \label{con:warehouseflow_1} \\
&\sum_{(h,j) \in \delta^+(h)}x_{hj}^{k} = \sum_{(j,h) \in \delta^-(h)}x_{jh}^{k} \quad \forall h \in H, k \in K \label{con:satelliteflow} \\
& \sum_{k \in K^1}\sum_{(h,j) \in \delta^+(h)}x_{hj}^k \geq l_h \quad \forall h \in H \label{satelliteOpening} \\
& \sum_{c \in C} D_c(w_{ch} + y_{ch}) \leq Q_h \quad \forall h \in H \label{con:satelliteCap}\\
& \sum_{k \in K^1}\sum_{(i,h) \in \delta^-(h)}z_{ih}^{k}-\sum_{k \in K^1}\sum_{(h,i) \in \delta^+(h)}z_{hi}^{k} = \sum_{c \in C}D_c(w_{ch}+y_{ch}) \quad \forall h \in H  \label{con:satelliteLoad}\\
& \sum_{k \in K^1}\sum_{(i,0) \in \delta^-(0)}z_{i0}^{k}-\sum_{k \in K^1}\sum_{(0,i) \in \delta^+(0)}z_{0i}^{k} = - \sum_{h \in H}\sum_{c \in C}D_c(w_{ch}+y_{ch}) \label{con:initialLoad}\\
& \sum_{\omega \in \Omega}f_{k\omega} = \sum_{h \in H}\sum_{(h,j) \in \delta^+(h)} x_{hj}^k \quad \forall k \in K^1 \label{con:stopsCount}\\
& z_{h0}^{k} = 0 \quad \forall h \in H, k \in K^1 \label{con:setLoad} \\
&\text{Echelon 2:}\notag\\
& \sum_{k \in K^2}\sum_{(c,j) \in \delta^+(c)}x_{cj}^{k} = \sum_{k \in K^2}\sum_{(j,c) \in \delta^-(c)}x_{cj}^{k} \quad \forall c \in C \label{con:customersflow1} \\
& \sum_{k \in K^2}\sum_{(j,c) \in \delta^-(c)}x_{cj}^{k} = 1 - \sum_{h \in H_c}w_{ch} \quad \forall c \in C \label{con:customersflow2} \\
& w_{ch} + y_{ch} \leq l_h \quad \forall c \in C, h \in H_c \label{con:custTolocker} \\
& \sum_{h \in H}(w_{ch}+y_{ch}) = 1 \quad \forall c \in C \label{con:custOnesatellite}\\
& s_{c}^k \leq 1 - q_h^k +y_{ch} \quad \forall c \in C, h \in H_c, k \in K^2\label{con:r_y_connect1}\\
& s_{c}^k \leq 1 + q_h^k -y_{ch} \quad \forall c \in C, h \in H_c, k \in K^2\label{con:r_y_connect2}\\
& \sum_{k \in K}s_{c}^k = \sum_{h \in H_c}y_{ch} \quad \forall c \in C \label{con:onepercustomer}\\
& \sum_{(h,j) \in \delta^+(h)}x_{hj}^{k} = \sum_{(j,h) \in \delta^-(h)}x_{jh}^{k}  \quad \forall h \in H, k \in K^2\label{con:satelliteflow2} \\
& \sum_{(j,h) \in \delta^-(h)}x_{jh}^{k} \leq  q_h^k  \quad \forall h \in H, k \in K^2\label{con:satelliteflow2_1} \\
& q_h^k \leq l_h \quad \forall h \in H, k \in K^2 \label{con:vehicletoopen}\\
& \sum_{h \in H} q_h^k \leq 1 \quad \forall k \in K^2 \label{con:maxonesatellite}\\
& \sum_{(i,c) \in \delta^-(c)}z_{ic}^{k}-\sum_{(c,i) \in \delta^+(c)}z_{ci}^{k} = D_cs_{c}^k \quad \forall c \in C, k \in K^2 \label{con:customerLoad}\\
& \sum_{(c,i) \in \delta^-(c)}x_{ci}^{k}+\sum_{(i,h) \in \delta^+(h)}x_{ih}^{k} - s_{c}^k \leq 1 \quad \forall c \in C, h \in H_c, k \in K^2 \label{con:linkcustomerlocker}\\
&  \sum_{k \in K^2}\sum_{(h,c) \in \delta^-(c)}z_{hc}^k = \sum_{c \in C}D_cy_{ch} \quad \forall h \in H\label{con:demandsatellite}\\
& w_{ch} = 0 \quad \forall c \in C, h \notin H_c \label{con:custHub}\\ 
& z_{ch}^{k} = 0 \quad \forall c \in C, h \in H, k \in K^2 \label{con:finalLoad}\\
& z_{ij}^{k} \leq Q_kx_{ij}^{k} \quad \forall (i,j) \in A, k \in K \label{con:veichlesCap}\\
& x_{ij}^k \in \{0,1\} \quad \forall (i,j) \in A, k \in K \label{domx}\\
& y_{ch} \in \{0,1\} \quad \forall c \in C, h \in H \label{domy}\\
& w_{ch} \in \{0,1\} \quad \forall c \in C, h \in H \label{domw}\\
& q_{h}^k \in \{0,1\} \quad \forall k \in K^2, h \in H \label{domq}\\
& s_{c}^k \in \{0,1\} \quad \forall c \in C, k \in K^2 \label{doms}\\
& l_{h} \in \{0,1\} \quad \forall h \in H \label{doml}\\
& z_{ij}^k \geq 0 \quad \forall (i,j) \in A, k \in K \label{domz}
\end{alignat}
The objective function \eqref{con:obj} aims to minimize the total $CO_2$ emissions generated by the delivery network and involves three components. The first component captures emissions related to the frequency of vehicle stops at intermediate satellites in the first-echelon. This is modeled through a piecewise linear function, where each additional stop made by vehicle $k \in K^1$ within an interval $\omega \in \Omega$ incurs emissions equal to $e_{\omega}$. The second component calculates routing emissions from vehicles across both echelons based on their emission factors ($e_k$) and the Fuel Consumption Rate (FCR) for each arc $(i,j)$, as detailed in constraints \eqref{eq:FCR}. Following the framework proposed by \cite{xiao}, parameters $\rho_k^*$ and $\rho_k$ represent fuel consumption rates when vehicles are fully loaded or empty, respectively. The final component accounts for emissions resulting from customers traveling to satellites.
Constraints \eqref{con:warehouseflow} through \eqref{con:setLoad} govern first-echelon operations. Constraints \eqref{con:warehouseflow} and Constraints \eqref{con:warehouseflow_1} ensure that each vehicle in $K^1$ may depart from the warehouse (node $0$) via at most one outgoing arc and return to the same node. Constraints \eqref{con:satelliteflow} manage flow conservation at intermediate satellites, while constraints \eqref{satelliteOpening} dictate that a satellite is considered active ($l_h=1$) only if at least one vehicle from $K^1$ visits it. Constraints \eqref{con:satelliteCap} limit the total demand assigned to each satellite to its respective capacity $Q_h$. Load flow conservation at satellites and the warehouse is enforced through constraints \eqref{con:satelliteLoad} and \eqref{con:initialLoad}, respectively. Specifically, \eqref{con:satelliteLoad} balance inbound and outbound flows at each satellite to match assigned demand, whereas \eqref{con:initialLoad} ensures net outflow from the warehouse equals the total demand across all satellites. Constraints \eqref{con:stopsCount} define the intervals for emissions calculations based on the number of arcs departing from satellites for each vehicle in $K^1$. Lastly, constraints \eqref{con:setLoad} guarantee that vehicles return empty to the warehouse at the conclusion of their routes.

Constraints \eqref{con:customersflow1} to \eqref{con:veichlesCap} control second-echelon operations. Constraints \eqref{con:customersflow1} and constraints \eqref{con:customersflow2} specify that a customer is served by a vehicle from $K^2$ only if they do not pick up their package from any satellite (i.e., $\sum_{h \in H_c}w_{ch}=0$). Constraints \eqref{con:custTolocker} and \eqref{con:custOnesatellite} ensure that each customer is assigned exclusively to one satellite and only if that satellite is activated. Constraints \eqref{con:r_y_connect1} and \eqref{con:r_y_connect2} link the assignment of second-echelon vehicles to satellites with customer assignments, ensuring consistency across these assignments. Constraints \eqref{con:onepercustomer} ensure that each customer is served by exactly one vehicle in $K^2$. Flow-balance constraints \eqref{con:satelliteflow2} and constraints \eqref{con:satelliteflow2_1} guarantee satellite visits by second-echelon vehicles only when explicitly assigned (i.e., $q_h^k=1$), an assignment further constrained by satellite activation and uniqueness via constraints \eqref{con:vehicletoopen} and \eqref{con:maxonesatellite}. Constraints \eqref{con:customerLoad} dictate load delivery to customers, ensuring that each customer's demand $D_c$ is met precisely by their assigned vehicle. Constraints \eqref{con:linkcustomerlocker} synchronize second-echelon routing with customer-satellite assignments, while constraints \eqref{con:demandsatellite} confirm that the total departing load from each satellite matches the combined demands of customers served through home delivery. Constraints \eqref{con:custHub} state that a customer $c$ can not travel to a satellite outside the available subset $H_c$ containing all the satellites at a distance equal or lower than $d_{max}^c$.
Constraints \eqref{con:finalLoad} ensure that vehicles in the second-echelon conclude their routes empty upon returning to their designated satellite. Finally, vehicle capacity constraints \eqref{con:veichlesCap} apply uniformly across both echelons, limiting load quantities carried on each arc to the respective vehicle capacities.

Variable domains are established by constraints \eqref{domx} through \eqref{domz}.
To further strengthen the proposed model, we incorporate the following valid inequalities.
The first inequality ensures that the total capacity of all active satellites is at least sufficient to cover the overall demand: \begin{alignat}{3} \sum_{h \in H} Q_h l_h \geq  \sum_{c \in C} D_c
\end{alignat}
The second set of inequalities enforces, for each satellite $h$, that the capacity of its assigned vehicles is large enough to meet the demand of the customers requiring home delivery and that are assigned to that satellite: \begin{alignat}{3} \sum_{k \in K^2}\sum_{h \in H} Q_kq_{h}^k \geq \sum_{c \in C} D_cy_{ch} \quad \forall h \in H \end{alignat}
Finally, we introduce a similar constraint focusing on customers who are not eligible to pick up orders, i.e., those customers $c$ with $d_{max}^c=0$. Let $C_0 \subseteq C$ denote this subset of customers, the corresponding valid inequality ensures that the vehicles assigned to satellites have enough capacity to serve all such customers: \begin{alignat}{3} \sum_{k \in K^2} Q_kq_{h}^k \geq \sum_{c \in C_0} \sum_{h \in H_c} D_cy_{ch} \end{alignat}

Having multiple vehicles with identical capacity and emission factor results in numerous symmetric solutions that unnecessarily increase the search space. To mitigate this issue, we introduce symmetry-breaking constraints to impose an explicit order among vehicles with the same attributes, specifically, we impose that lower index vehicles have to be used before the highest one. Let $K_\theta^\epsilon \subseteq K = \{k \in K: Q_k = \theta, e_k = \epsilon\}$ be the subset of vehicles having capacity $\theta$ and emission factor $\epsilon$, the symmetry-breaking constraints can be written in the form: 
\begin{alignat}{3}
    \sum_{h \in H} q_{h}^k \leq \sum_{h \in H} q_{h}^{k-1} \quad \forall k \in K_\theta^\epsilon\setminus\{0\}
\end{alignat}
In this way, vehicle $k$ is not assigned to any satellite if vehicle $k-1$ is not used.

\section{Design of Experiments} \label{sec:solution}
In this section, we detail the experimental setup designed to evaluate the 2E-LRP-ECB. First, in Section \ref{sec:scenarios} we describe four representative scenarios, each capturing a distinct configuration of objectives and fleet composition. Then, in Section \ref{sec:multi} we present the multi-objective analysis implemented to investigate trade-offs between different objectives.

\subsection{Scenario-based Analysis} \label{sec:scenarios}
To thoroughly investigate the environmental impact of two-echelon location-routing configurations, we have solved four different scenarios of the proposed 2E-LRP-ECB model, each reflecting a distinct strategic operational policy. These scenarios differ based on their primary optimization objective and the capacity of the green fleet available in the second-echelon. Specifically, the scenarios considered are:

\begin{enumerate}
\item[EHC:] \textbf{Emission minimization with High-Capacity green vehicles.} 
In this scenario, overall emissions are minimized while equipping the second echelon with higher-capacity green vehicles. 
The goal is to assess how larger vehicle capacity influences emission levels, satellite activation, and customer pick-ups when environmental sustainability is the dominant concern. 
By prioritizing low emissions, this scenario indirectly emphasizes social sustainability (e.g., improved air quality) and may also reduce economic costs linked to fuel consumption.
\\
\textbf{Objective function:} As formulated in the model in equation \eqref{con:obj}.

\item[ELC:] \textbf{Emission minimization with Low-Capacity green vehicles}: This scenario also targets overall emissions minimization but restricts the second-echelon to lower-capacity green vehicles. 
While additional trips by electric vehicles do not increase emissions, their reduced capacity may limit the system’s ability to serve all demand exclusively with green vehicles. 
As a consequence, it may be necessary to assign some deliveries to combustion vehicles, which reintroduces emissions into the system. 
The motivation is therefore to analyze the environmental and operational trade-offs that arise when green vehicle capacity is constrained, particularly the balance between the usage of green and combustion fleets.
\\
\textbf{Objective function:} As formulated in the model in equation \eqref{con:obj}.

\item[TD:] \textbf{Total Distance minimization}: Here, the objective is to minimize the total traveled distances considering both company operations and customer travels to pick-up points. The rationale behind this scenario is to evaluate a situation where the delivery company aims to enhance overall customer convenience and operational efficiency by collectively reducing both the company’s and customers’ travel distances.
\\
\textbf{Objective function:} The function in equation~\eqref{con:obj} is replaced with
\begin{equation}
    \sum_{k \in K}\sum_{(i,j) \in A} d_{ij}x_{ij}^k + \sum_{c \in C}\sum_{h \in H} d_{ch}w_{ch}, \label{con:objDist}
\end{equation}
where the goal is to minimize both company and customer travel distances.

\item[CD:] \textbf{Company Distance minimization}: This scenario solely minimizes distances traveled by the company's vehicles, excluding customers' trips. The motivation is to assess the effectiveness of focusing purely on company-level operational efficiency without considering customer-side trips.
\\
\textbf{Objective function:} The function in equation \eqref{con:obj} is replaced with
\begin{equation}
    \sum_{k \in K}\sum_{(i,j) \in A} d_{ij}x_{ij}^k, \label{con:objCompDist}
\end{equation}
with the aim of minimizing company travel distances only.
\end{enumerate}

For the emission-oriented scenarios (EHC and ELC), the objective function remains unchanged from the base formulation in Section \ref{sec:model}, with the only distinction being the available capacity of the green fleet. 
In contrast, the distance-oriented scenarios (TD and CD) adopt modified formulations: the former replaces the objective function with equation \eqref{con:objDist}, capturing both company and customer trips, while the latter employs equation \eqref{con:objCompDist}, focusing solely on company travel. 
Across all scenarios, the set of constraints from \eqref{con:warehouseflow} to \eqref{domz} is preserved without modification.

After solving each scenario exactly, we perform a comparative analysis to highlight key differences and provide managerial insights. Subsequently, a multi-objective analysis considering emissions and distances (as described in Section \ref{sec:multi})  is conducted to clearly illustrate trade-offs between environmental sustainability and operational efficiency.
The results of these analyses will allow decision-makers to better understand the implications of strategic choices related to fleet composition, satellite activation, and the incorporation of customer behaviors into the last-mile delivery system.

\subsection{Multi-Objective approach} \label{sec:multi}
In multi-objective optimization, decision-makers often face trade-offs between conflicting objectives. In our setting, we investigate the trade-off between minimizing total emissions and minimizing total traveled distances within a capacitated two-echelon location-routing framework. These two objectives often do not go in the same optimization direction due to the possibility for customers to pick up their orders with zero-emissions modes of transportation and for green second-echelon vehicles. These two options only increase the traveled distance without adding emissions to the solution.
To rigorously explore this trade-off, we implement the $\epsilon$-constraint method \cite{MesquitaCunha2023} a well-established scalarization technique for generating Pareto-optimal solutions in multi-objective problems. 

The $\epsilon$-constraint method consists of optimizing one objective while transforming the others into additional constraints bounded by predefined thresholds ($\epsilon$ values). Specifically, given a multi-objective optimization problem with several objective functions $f_1(x), f_2(x), \dots, f_n(x)$, the method selects one of the functions, for instance, $f_1(x)$, to minimize:

\begin{align*}
\min & \quad f_1(x) \\
\text{s.t.} & \quad f_j(x) \le \epsilon_j, \quad j = 2,\dots,n \\
& \quad x \in X
\end{align*}

where $X$ represents the set of feasible solutions and $\epsilon_j$ are parameters systematically varied to explore different regions of the Pareto frontier.

In this study, we apply the $\epsilon$-constraint method to quantify the trade-offs between emissions and distances. Initially, we independently minimize each objective to obtain their respective ideal solutions, thereby identifying the range of values each objective can assume. We then progressively vary the value of $\epsilon$, constraining the secondary objective, and solve the resulting mixed integer linear programming (MILP) problems exactly, generating a comprehensive set of Pareto-optimal solutions.

The resulting Pareto frontier provides critical insights into the environmental and operational implications of emphasizing either emissions reduction or distance minimization. Each Pareto-optimal point represents a solution in which improving one objective necessarily leads to a compromise in the other, guiding decision-makers toward balanced and informed operational choices.
Once the frontier have been drawn, we compute the knee-point. The knee-point is determined with a geometric approach by calculating the point on the frontier with the maximum perpendicular distance from a straight line connecting the extreme solutions (in our case, the points with minimum emissions and minimum distances).
Specifically, let $(x_i,y_i)_{i=1}^N$ be the Pareto‐optimal points, ordered so that
$((x_1,y_1)$ and $(x_N,y_N)$ are the extremes (minimum emissions and minimum distances).  
For each intermediate point $i=\{2,\dots,N-1\}$, compute its perpendicular distance $d_i$ to the line joining the extremes as :
\begin{equation}
      d_i = \frac{\bigl|(x_N - x_1)\,(y_1 - y_i)\;-\;(x_1 - x_i)\,(y_N - y_1)\bigr|}
             {\sqrt{(x_N - x_1)^2 + (y_N - y_1)^2}}.
\end{equation}

The knee‐point index $k$ is then
$
  k = \arg\max_{i=2,\ldots,N-1} d_i,
$
and the knee‐point coordinates are $(x_k,y_k)$.

\section{Benchmark Instances based on a Realistic Case Study}\label{sec:casestudy}
\begin{figure}[htbp!]
  \centering
  \scalebox{0.8}{
  \includegraphics[width=0.6\textwidth]{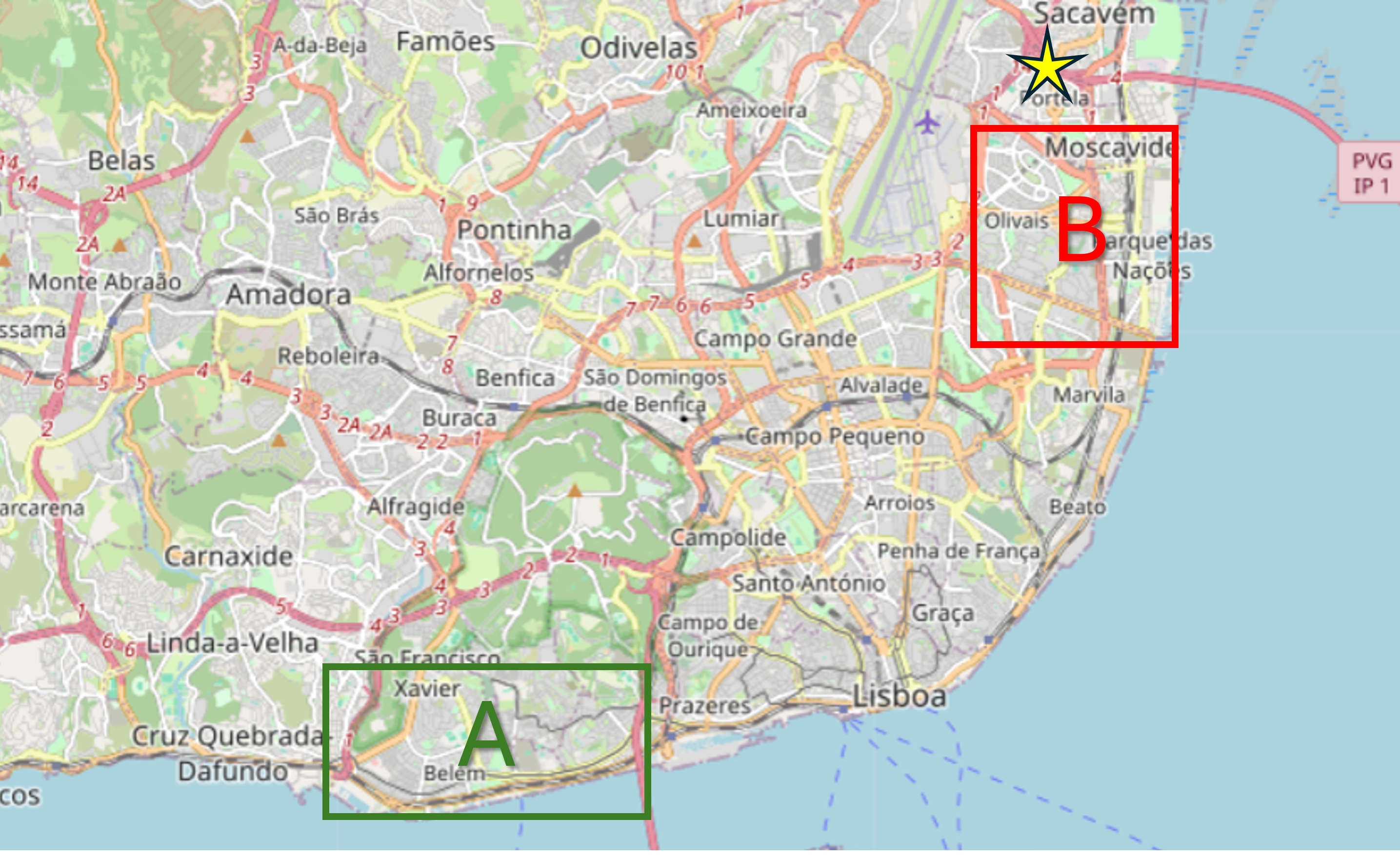}}
  \caption{Map of Lisbon}
  \label{fig:map}
\end{figure}

\begin{figure}[htbp!]
  \centering
  \begin{minipage}[b]{0.48\linewidth}
    \centering
     \scalebox{0.9}{
    \includegraphics[width=\linewidth]{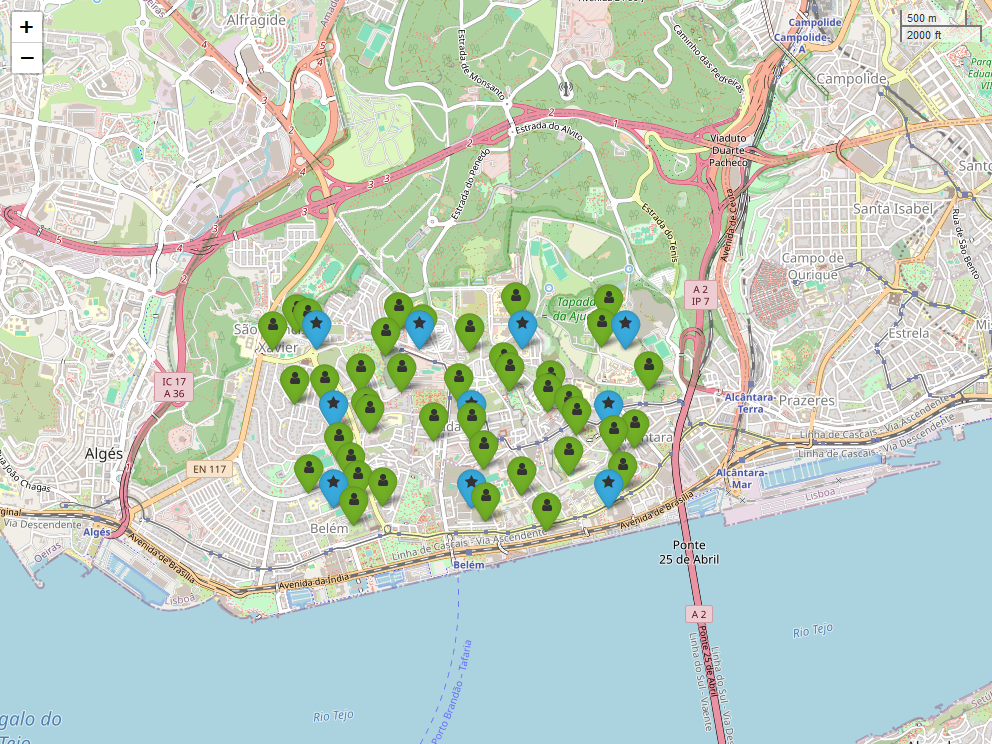}}
    \caption{Potential satellites and illustrative customer distribution in Area A. }
    \label{fig:1300}
  \end{minipage}
  \hfill
  \begin{minipage}[b]{0.48\linewidth}
    \centering
     \scalebox{0.9}{
    \includegraphics[width=\linewidth]{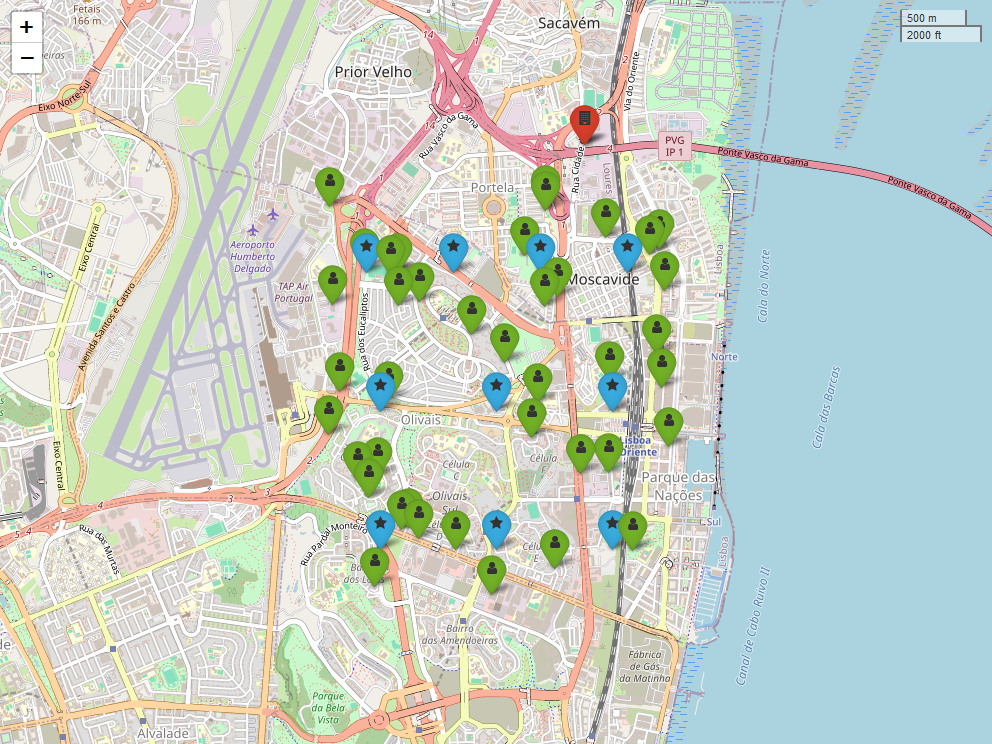}}
    \caption{Potential satellites and illustrative customer distribution in Area B.}
    \label{fig:1800}
  \end{minipage}
\end{figure}

To test and validate our approach to the 2E-LRP-ECB, we generated benchmark instances based on a realistic case study of Lisbon, Portugal. For the first-echelon, we considered a warehouse located near the east side of the city (Oriente), serving as the main hub for city operations. From this warehouse, orders are dispatched to multiple satellites distributed across the city. To balance computational complexity with robustness and generalizability of the results, we selected two distinct areas for detailed analysis:
\begin{itemize}
    \item Area A, located farther form the warehouse and covering the west area of the city (Belém).
    \item Area B, situated in the east area of the city, closer to the warehouse.
\end{itemize}
This choice enables us to explore and compare the impact of geographic proximity and logistical complexity on emissions, operational costs, and customer behaviors. Figure \ref{fig:map} shows the city of Lisbon considered in this study, with the warehouse marked by a star and the two areas highlighted with rectangular shapes. Figures \ref{fig:1300} and \ref{fig:1800} provide details of the satellites (blue symbols) and an example of customer distribution in Areas A and B, respectively.

\hfill

Considering the number of customers, their locations, and demand, we generated 80 instances. Table \ref{tab:instances} presents their characteristics. The columns list the number of customers ($|C|$), the number of available satellites ($|H|$), and the different vehicle sets, while the analyzed areas are indicated in the “Area” column. To mimic real-world demand fluctuations, we varied customer-related inputs (e.g., locations and demands), which influence satellite activation. In practice, satellites such as pickup stores and lockers remain available at fixed locations; therefore, we kept their geographic coordinates constant. However, the actual activation of each satellite, as well as the customers assigned to them, depends on demand, creating a dyamic-satellite configuration. Based on this, we study how strategic choices regarding satellite utilization affect both routing costs and customer service levels.
\begin{table}[htbp]
\centering
\scriptsize
\caption{Instance parameters by problem size}
\begin{tabular}{lrrrr}
\toprule
Size & $\lvert C\rvert$ & $\lvert H\rvert$ &$Area$& \#\textbf{Instances} \\
\midrule
Small  & \{40\}& 10 & \{A,B\} & 10 \\
Medium & \{50,75\} & 10 &\{A,B\} & 10 \\
Large  &  \{100,150\} & 10 &\{A,B\}& 5 \\
\bottomrule
\end{tabular}
\label{tab:instances}
\end{table}
For all instance sets, customer locations were randomly generated within the study area, which was divided into $|H|$ zones, each containing one satellite. Each satellite has the capacity to accommodate around 15 packages on average. Customers’ emission profiles were generated using five values of $d_{max}$ (0 m, 500 m, 1000 m, 1500 m, and 2000 m) and two values of $d_{green}$ (0 m, 500 m). These combine into nine meaningful [$d_{green}$,$d_{max}$] configurations, with $d_{green} \leq d_{max}$. When both values are zero, it indicates that the customer is unwilling to travel to pick up an order.
Customer demand was drawn from four package size categories (extra-small, small, medium, and large) as summarized in Table \ref{tab:packages}. We assumed default home delivery for extra-small and large categories: the former because such parcels fit directly into mailboxes, and the latter because their weight precludes customer pickup.
\begin{table}[htbp!]
\centering
\caption{order size categories and their frequencies}
\begin{tabular}{l r r}
\toprule
Type & Size (up to) & Distribution \\
\midrule
Extra-small & 0.35 kg     & 30\% \\
Small       & 2 kg       & 40\% \\
Medium      & 5 kg  & 25\% \\
Large       & 30 kg   &  5\% \\
\bottomrule
\end{tabular}
\label{tab:packages}
\end{table}

Other important parameters are those related to vehicles (types, capacities, and emissions), which are summarized in Table 8—specifically, vehicle capacity ($Q_k$), emission factors ($e_k$), and fuel consumption rates ($FCR_K$). We consider three vehicle types common to all instances, differing in both capacity and emissions. First-echelon trucks are denoted as $K^1$, while $K^2$ and $K^2_0$ represent, respectively, combustion vans and zero-emission vehicles. $K^2_0$ includes two different capacities: they can accommodate on average 25 packages in the ELC scenario and 50 for all others. For customers traveling with combustion vehicles, the emission factor is set at 0.15 kg/km. All $CO_2$ emission data (for commercial and private vehicles) were obtained from the EMEP/EEA Guidebook \cite{EMEP2024}.
Drawing on practical experience, each satellite stop is assumed to last five minutes. Additional emissions per stop are modeled using a piecewise function \cite{EMISIA_COPERT5_2025}, which accounts for variability in stopping durations and conditions: 0.1 kg $CO_2$/stop for fewer than two stops, 0.15 kg $CO_2$/stop for two to four stops, and 0.3 kg $CO_2$/stop for five or more stops.

\begin{table}[htbp!]
  \centering
  \caption{Vehicles characteristics}
    \begin{tabular}{lrr}
    \toprule
   Vehicle  &\multicolumn{1}{r}{$e_k$ ($ kg/ Km$)} & \multicolumn{1}{r}{$\rho^*_k-\rho_k$ (l/Km)} \\
    \midrule
    $K^1$ & 0.38  & 0.5 \\
    $K^2$ & 0.30 & 0.2 \\
    $K^2_0$ & 0 & 0 \\
    \bottomrule
    \end{tabular}%
  \label{tab:vehicles}%
\end{table}%

\section{Computational Results} \label{sec: settings}
In this section, we present the computational results obtained by applying the methodology to the case-study inspired instances described in Section \ref{sec:casestudy}. Section \ref{sec:envres} analyzes the distribution of total $CO_2$ emissions and traveled distances across the two-echelon vehicles and customers picking up their orders, providing insights into the changes that occur when shifting from an emission-minimization objective (EHC scenario) to a distance-minimization one (TD scenario). Section \ref{sec:sat} examines the spatial and operational activation of satellites, while Section \ref{sec:customers} reports statistics on customer behavior. Finally, Section \ref{sec:multires} presents a set of Pareto frontiers obtained by jointly considering emissions and distance minimization.
All tests were conducted on a machine running Windows 10 Pro, equipped with an Intel Xeon E5-2660 TD CPU (20 cores) and 64 GB of RAM. The MILP model was implemented in Java 21 and solved using Gurobi 11.0.3. The average MIP gap and computational times for each instance size and problem scenario are reported in Table \ref{tab:gap}. The Gurobi time limit was set to 3600 s for small and medium-size instances and 18000 s for large ones. All small- and medium-size instances were solved to optimality across all scenarios, while among the large instances, only those with 150 customers failed to reach optimality in any scenario.

\begin{table}[htbp!]
\centering
\scriptsize
\caption{Average computational times and gaps for instance sizes and problem scenarios}
\begin{tabular}{lrr|rr|rr|rr}
\toprule
& \multicolumn{2}{c|}{EHC} & \multicolumn{2}{c|}{ELC} & \multicolumn{2}{c|}{TD} & \multicolumn{2}{c}{CD} \\
$|C|$ & Time (s) & Gap (\%) & Time (s) & Gap (\%) & Time (s) & Gap (\%) & Time (s) & Gap (\%) \\
\midrule
40  &  31.7   & 0.0 &   78.1   & 0.0 &  122.3   & 0.0 &  129.5   & 0.0 \\
50  &  90.2   & 0.0 &  621.2   & 0.0 &  601.5   & 0.0 & 2010.4   & 0.0 \\
75  & 663.1   & 0.0 & 1285.1   & 0.0 & 4575.3   & 0.0 & 4940.6   & 0.0 \\
100 & 4738.7  & 0.0 & 6926.5   & 0.0 & 15743.1  & 2.9 & 16326.9  & 3.1 \\
150 & 18000.0 & 7.6 & 18000.0  & 3.1 & 18000.0  & 9.8 & 18000.0  & 9.5 \\
\bottomrule
\end{tabular}
\label{tab:gap}
\end{table}

\subsection{Environmental Impact Analysis} \label{sec:envres}
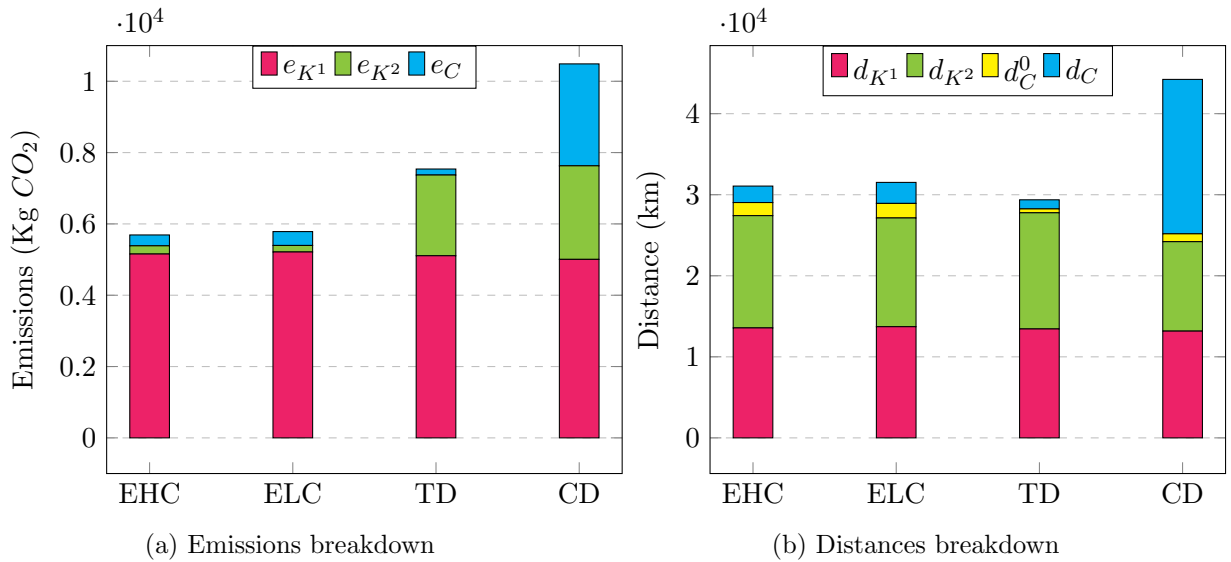
\begin{figure}[htbp]
\centering
\begin{minipage}[b]{0.48\textwidth}
  \centering
\begin{tikzpicture}
\begin{axis}[
    ybar stacked,
    bar width=15pt,
    enlargelimits=0.1,
    ylabel={Emissions (Kg $CO_2$)},
    symbolic x coords={EHC, ELC, TD, CD},
    xtick=data,
    legend style={at={(0.5,1)}, anchor=north, legend columns=-1},
    ymajorgrids=true,
    grid style=dashed,
    ymin=0,
    ymax=10000
]

\addplot[fill=WildStrawberry] coordinates {
    (EHC,5155.821)
    (ELC,5212.6196)
    (TD,5107.2152)
    (CD,5005.8046)
};

\addplot[fill=LimeGreen] coordinates {
    (EHC,229.65)
    (ELC,183.177)
    (TD,2263.563)
    (CD,2623.227)
};

\addplot[fill=ProcessBlue] coordinates {
    (EHC,304.794)
    (ELC,388.953)
    (TD,168.4425)
    (CD,2858.073)
};
\legend{$e_{K^1}$, $e_{K^2}$, $e_{C}$}
\end{axis}
\end{tikzpicture}
\subcaption{Emissions breakdown}
\label{fig:emissions}
\end{minipage}
\hfill
\begin{minipage}[b]{0.48\textwidth}
  \centering
  \begin{tikzpicture}
  \begin{axis}[
      ybar stacked,
      bar width=15pt,
      enlargelimits=0.1,
      ylabel={Distance (km)},
      symbolic x coords={EHC, ELC, TD, CD},
      xtick=data,
      ymin=0, ymax=44000,
      legend style={at={(0.5,1)}, anchor=north, legend columns=-1},
      ymajorgrids=true,
      grid style=dashed
  ]
  \addplot[fill=WildStrawberry] coordinates {(EHC,13567.95) (ELC,13717.42) (TD,13440.04) (CD,13173.17)};
  \addplot[fill=LimeGreen] coordinates {(EHC,13852.04) (ELC,13415.79) (TD,14333.47) (CD,11035.32)};
  \addplot[fill=Yellow] coordinates {(EHC,1616.02) (ELC,1796.16) (TD,477.98) (CD,968.45)};
  \addplot[fill=Cyan] coordinates {(EHC,2031.96) (ELC,2593.02) (TD,1122.95) (CD,19053.82)};
  \legend{$d_{K^1}$, $d_{K^2}$, $d_{C}^0$, $d_{C}$}
  \end{axis}
  \end{tikzpicture}
  \subcaption{Distances breakdown}
\label{fig:dist}
\end{minipage}
\caption{Emissions and distances breakdown by scenario (EHC: Emissions minimization with High-Capacity; ELC: Emissions minimization with Low Capacity; TD: Total Distances; CD: Company Distances)}
\label{fig:em-dist-breakdown}
\end{figure}

This section provides a comprehensive assessment of the environmental implications of the proposed model. 
Figure \ref{fig:em-dist-breakdown} illustrates the breakdown of total emissions and traveled distances for different problem scenarios. These values are obtained averaging the results among all instances.
In Figure \ref{fig:emissions}, the column representing total emissions is divided into three main sources: first-echelon vehicles ($e_{K^1}$), second-echelon vehicles ($e_{K^2}$) and customers traveling ($e_{C}$). 
Figure \ref{fig:dist} breaks the total traveled distances into four main sources: first-echelon vehicles ($d_{K^1}$), second-echelon vehicles ($d_{K^2}$), customers picking up the orders without any emissions ($d_{C}^0$) and with emissions ($d_{C}$).

In Figure \ref{fig:emissions}, it is evident that first-echelon vehicles consistently account for the largest share of total emissions across all scenarios. When emissions minimization is prioritized (in both EHC and ELC), second-echelon emissions are mainly attributable to customers. This outcome arises because it is more beneficial to rely on customer trips to satellites rather than employing additional commercial vehicles. In these scenarios, zero-emission vehicles are fully utilized in the second echelon, and whenever their capacity is insufficient, customers travel to satellites to shorten the routes of emission-based company vans.
In contrast, when the objective shifts to minimizing total traveled distances (TD), emissions increase considerably. This occurs because the optimization model focuses on routing efficiency without fully internalizing the emissions implications of vehicle and customers heterogeneity, resulting in greater reliance on second-echelon vehicles with non-zero emission profiles. The extreme case is represented by scenario CD, where distance minimization excludes customer travels from the objective function, causing an exponential rise in emissions.

Figure \ref{fig:dist} complements this analysis by showing how distances traveled in the first and second-echelon remain relatively balanced across scenarios. Moreover, distances traveled by customers using combustion vehicles systematically exceed those using zero-emission ones due to larger permitted travel radius ($d_{max}$). This can be linked to the model tending to reduce the number of active satellites to contain first-echelon emissions, often leaving customers without nearby collection points, thereby increasing average customer travel distances when pick-up is encouraged.
Both figures underscore the fundamental importance of including customers in the analysis to achieve a globally efficient system that avoids increases in emissions and traffic congestion, as shown in CD. 
Additionally, the results show that transitioning from EHC to ELC does not significantly modify system performance: the reduction in green-vehicle capacity in ELC is largely compensated by increased customer trips. Consequently, for clarity and focus, subsequent analyses will concentrate on directly comparing EHC and TD, which represent the clearest contrast between emissions-driven and distance-driven strategies.

\begin{figure}[htbp]
  \centering
  \begin{subfigure}{0.48\textwidth}
    \centering
    \includegraphics[width=\linewidth]{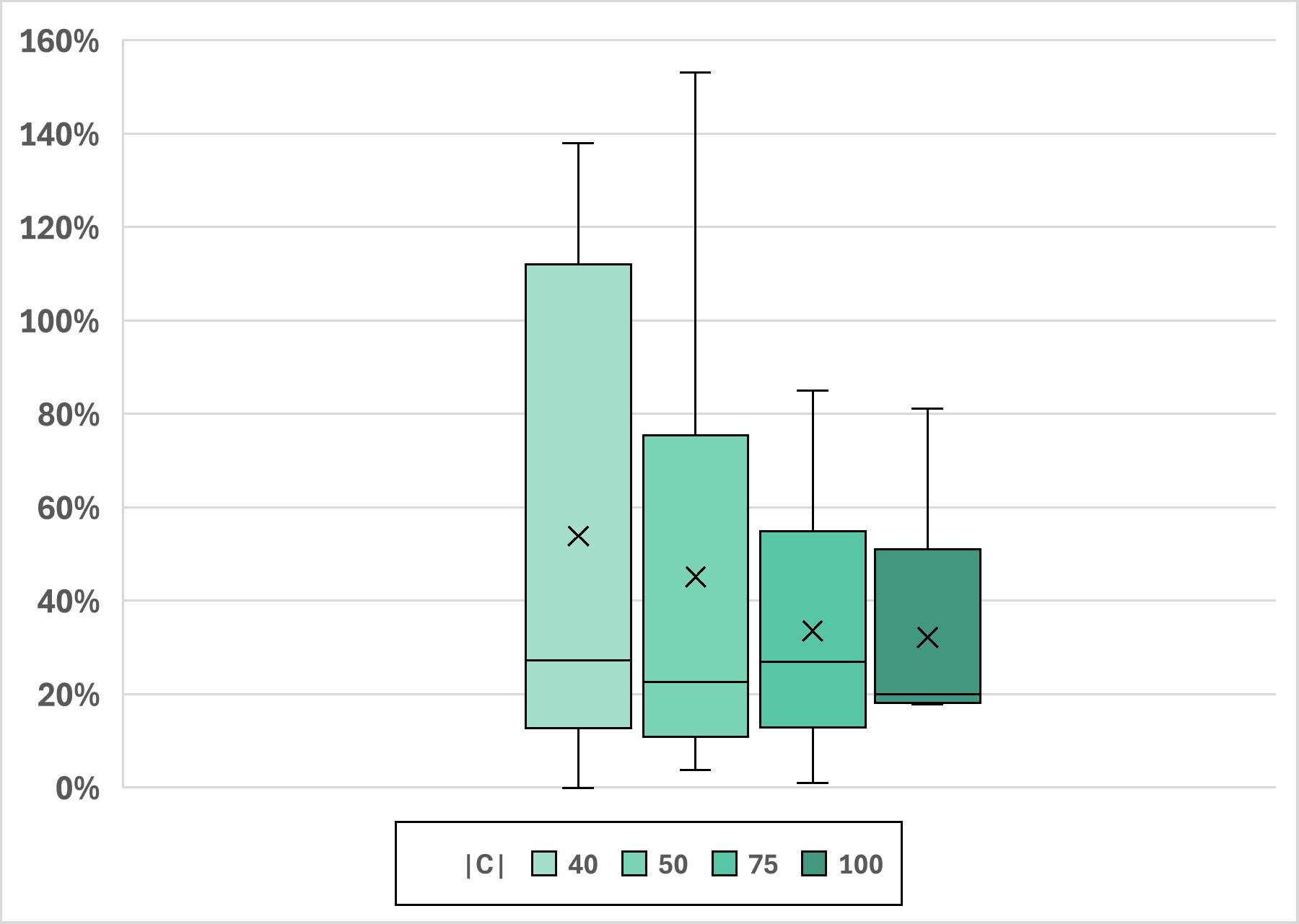} 
    \caption{Emissions increase}
    \label{fig:boxem}
  \end{subfigure}\hfill
  \begin{subfigure}{0.48\textwidth}
    \centering
    \includegraphics[width=\linewidth]{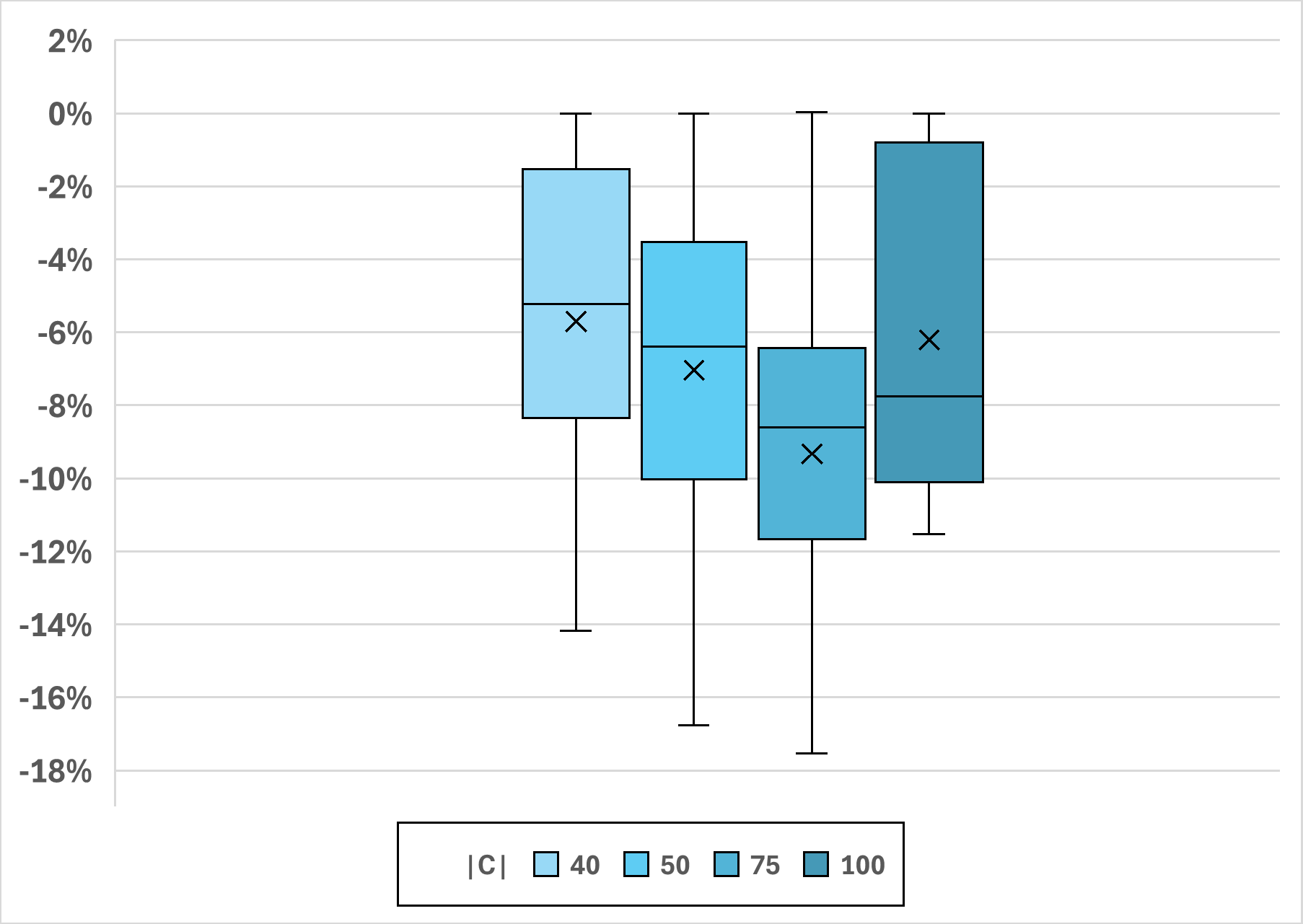} 
    \caption{Distances decrease}
    \label{fig:boxdist}
  \end{subfigure}
  \caption{Variation in emissions and traveled distances when moving from EHC (minimization of emissions) to TD (minimization of distances) for different number of customers ($|C|$)}
  \label{fig:boxplots}
\end{figure}

To better understand the increase in emissions and the decrease in traveled distances moving from EHC to TD, Figure \ref{fig:boxplots} shows the percentage variations for different instances size. Emissions increase and distances decrease are shown in Figure \ref{fig:boxem} and Figure \ref{fig:boxdist}, respectively.

Figure \ref{fig:boxem} presents a clear increase in emissions. For example, the average increase is about 20$\%$ in 40 customer instances and 17$\%$ in instances with 100 customers. The dispersion is greater in smaller instances  and progressively narrows as size grows since an increase in customers number promote the use of combustion vehicles. Figure \ref{fig:boxdist} displays the decrease in distances with median variations lying between -6$\% $ and -10$\% $ across all instances sizes.
These patterns underline a non-symmetrical trade-off: switching to a distance-oriented objective yields only modest savings in traveled kilometers, while it can provoke substantial increases in emissions, particularly in medium-size instances where the use of zero-emission vehicles is fundamental. From a managerial standpoint, this suggests that (i) pure distance minimization is difficult to justify, given that the average decrease in distances is always below 10\% whereas emissions may almost double and (ii) for larger-size instances, where variability is lower, $\epsilon$ -constraint methods can be employed to fine-tune the balance without incurring excessive environmental costs.

\subsection{Activation of Satellites}\label{sec:sat} 
In this section we discuss the activation of satellites in different resolutions. 
At first, Table \ref{tab:sat} gives an overview on someg satellites statistics: average number of active satellites with ($Active$) and without ($Active_{PU}$) a second-echelon route assigned and average number of customers per satellite ($AvgC$). Note that, active satellites without a second-echelon route assigned to them are destined only for order pick-up.

\begin{table}[htbp!]
\centering
\scriptsize
\caption{Satellites statistics for instance sizes and problem scenarios}
\begin{tabular}{lrrr|rrr|rrr|rrr}
\toprule
& \multicolumn{3}{c|}{EHC} & \multicolumn{3}{c|}{ELC} & \multicolumn{3}{c|}{TD} & \multicolumn{3}{c}{CD} \\
$|C|$ & Active & Active$_{\mathrm{PU}}$ & AvgC & Active & Active$_{\mathrm{PU}}$ & AvgC & Active & Active$_{\mathrm{PU}}$ & AvgC & Active & Active$_{\mathrm{PU}}$ & AvgC \\
\midrule
40  &  2.8 & 1.2 & 15.2 &  2.8 & 0.9 & 15.8 &  2.8 & 1.2 & 15.3 &  3.0 & 1.2 & 14.7 \\
50  &  3.8 & 2.2 & 14.2 &  3.4 & 1.4 & 16.3 &  3.9 & 2.3 & 14.1 &  3.5 & 1.5 & 15.3 \\
75  &  6.5 & 4.9 & 11.7 &  6.2 & 4.2 & 12.3 &  6.5 & 4.9 & 11.8 &  5.2 & 3.2 & 15.3 \\
100 &  8.1 & 6.6 & 12.4 &  8.5 & 6.5 & 11.8 &  8.0 & 6.3 & 12.5 &  6.3 & 4.3 & 16.1 \\
150 &  9.8 & 8.0 & 15.3 & 10.0 & 8.0 & 15.0 & 10.0 & 8.0 & 15.0 &  9.5 & 7.5 & 15.8 \\
\bottomrule
\end{tabular}
\label{tab:sat}
\end{table}

\begin{figure}[htbp!]
  \centering
  \resizebox{\linewidth}{!}{%
    \begin{tabular}{@{}cc@{}}
      \begin{subfigure}{0.48\linewidth}
        \centering
        \includegraphics[width=\linewidth]{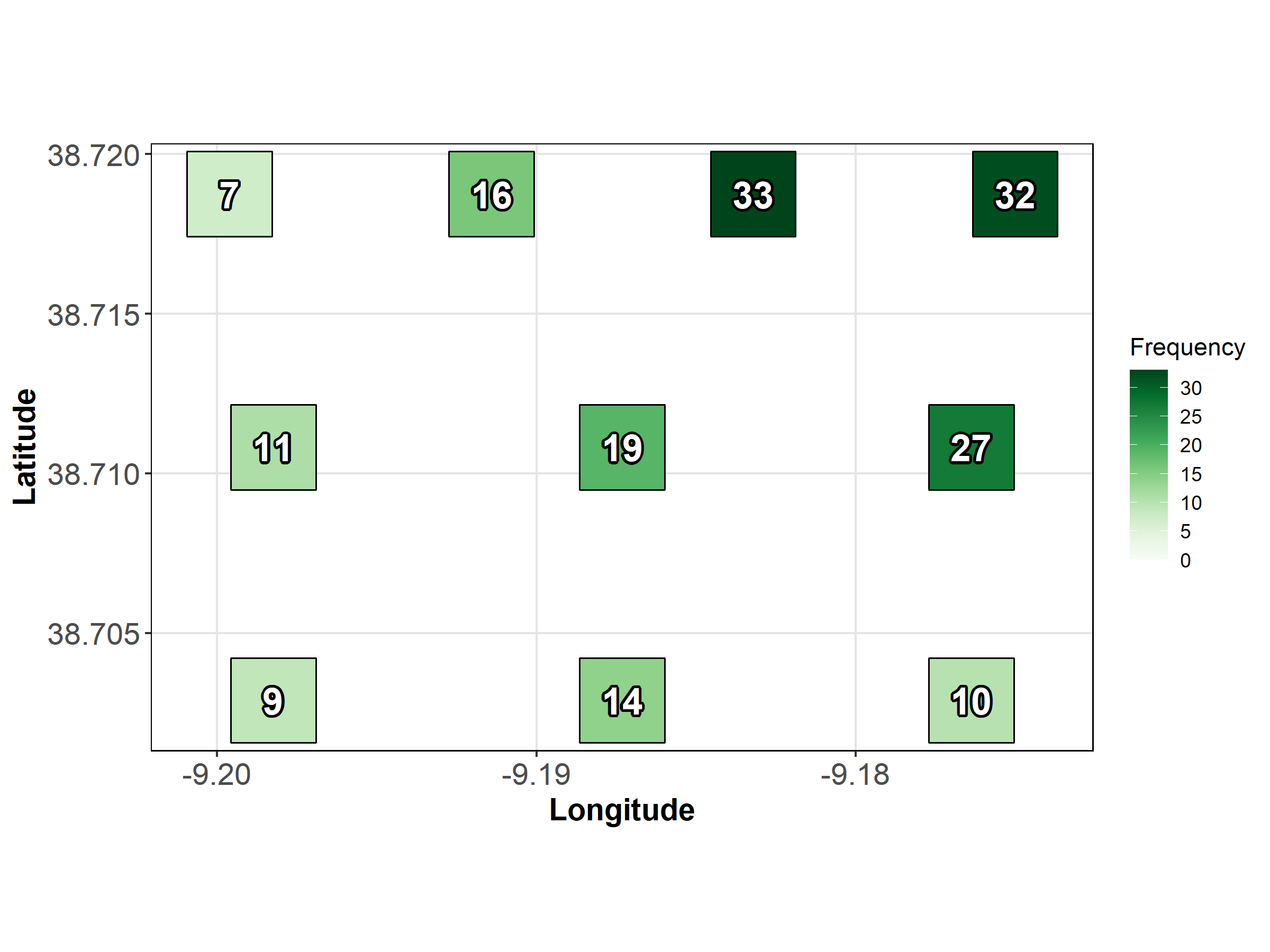}
           \caption{Area A - EHC}
      \label{fig:13300}
      \end{subfigure} &
      \begin{subfigure}{0.48\linewidth}
        \centering
        \includegraphics[width=\linewidth]{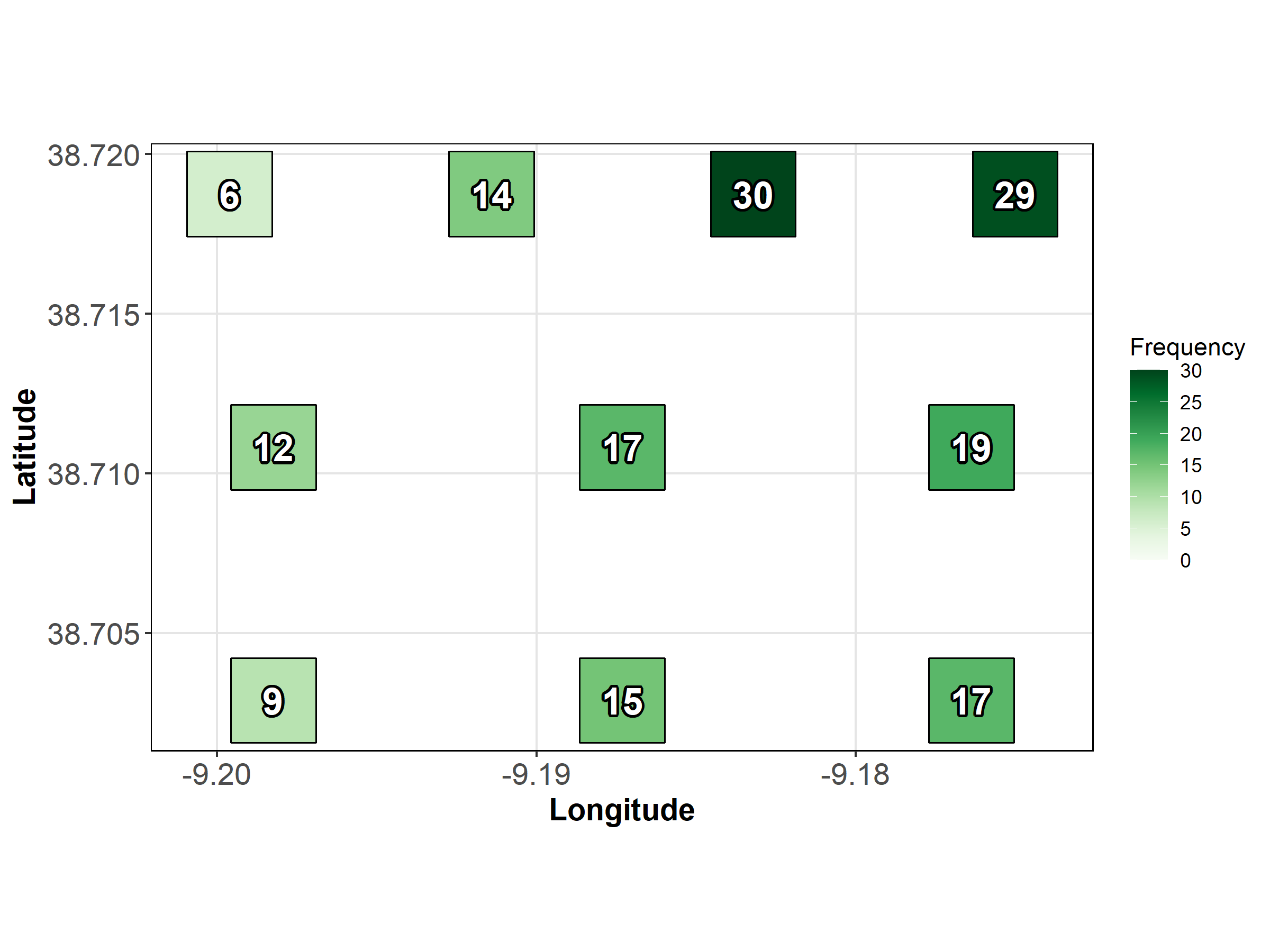}
      \caption{Area A - TD}
      \label{fig:13001}
      \end{subfigure} \\
      \begin{subfigure}{0.48\linewidth}
        \centering
        \includegraphics[width=\linewidth]{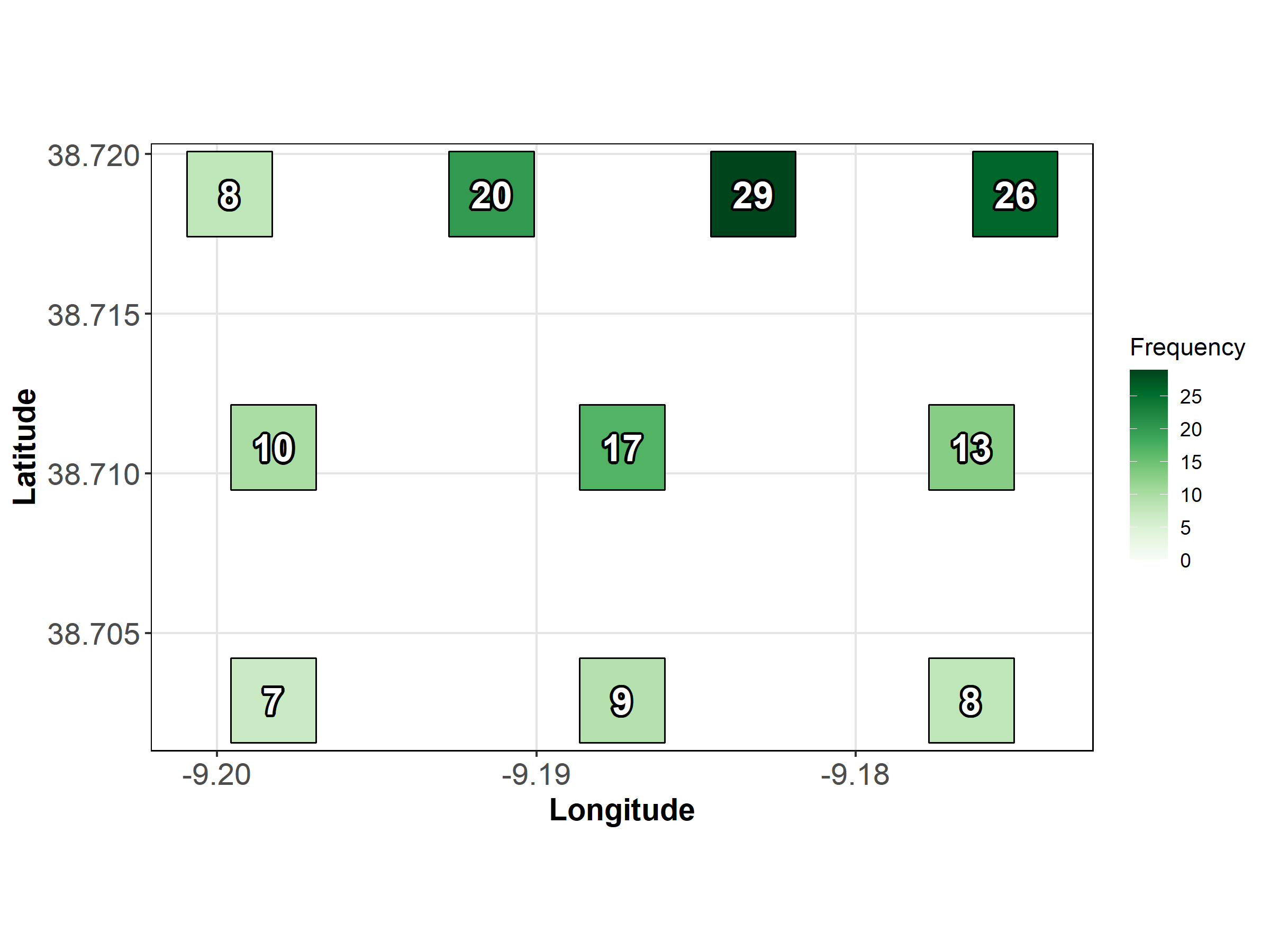}
              \caption{Area B - EHC}
      \label{fig:18000}
      \end{subfigure} &
      \begin{subfigure}{0.48\linewidth}
        \centering
        \includegraphics[width=\linewidth]{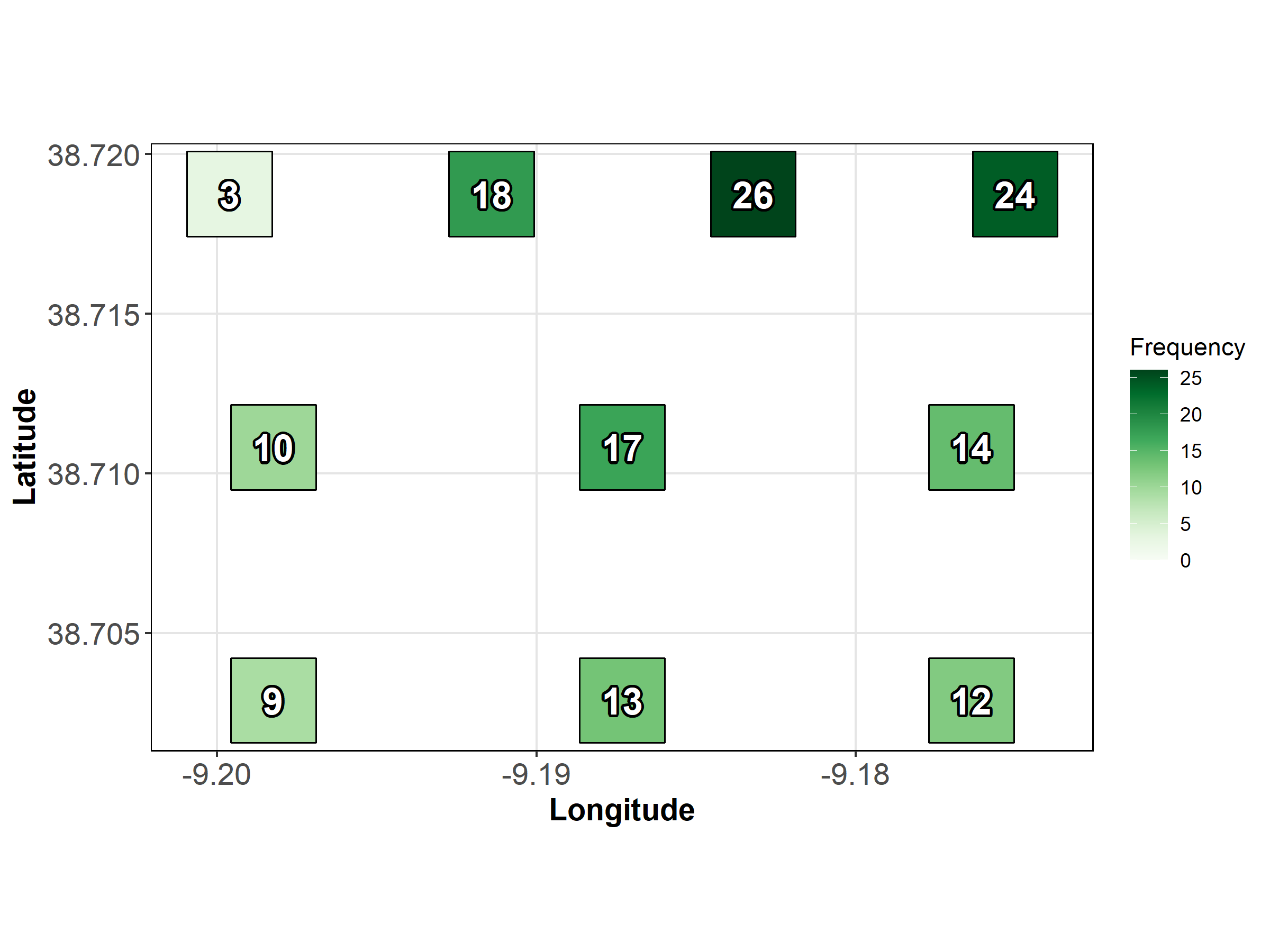}
           \caption{Area B - TD}
      \label{fig:18001}
      \end{subfigure}
    \end{tabular}%
  }
  \caption{Frequency of active satellites across all instances}
  \label{fig:heatmaps}
\end{figure}

Figure \ref{fig:heatmaps} presents a visual representation of the frequency of utilizing different satellites over a total of 35 instances for each area. Instances with $|C| = 150$ have been excluded, since in these larger cases the total demand exceeds the capacity of only a subset of satellites, making it necessary to activate all available satellites to satisfy customer needs. Figure \ref{fig:13300} and Figure \ref{fig:13001} show area A in EHC and TD, respectively, while Figure \ref{fig:18000} and Figure \ref{fig:18001} analyze area B, reflecting our decision to concentrate on EHC as the representative emissions-oriented strategy and TD as the distance-oriented one. The warehouse is located in the upper-right corner of every panel. 
This analysis highlights how daily activities depend on fluctuations in demand: the set of satellites activated changes from one day to another, reflecting a dynamic-satellite perspective where the network dynamically adapts to the spatial distribution of demand. This variability provides managerial insights on which satellites are consistently valuable and which serve only as occasional support.

Across all figures, the highest frequencies systematically belong to the satellites closest to the warehouse, confirming that the location of satellites tends to prioritize the reduction of first-echelon emissions. Under EHC (in both Figure \ref{fig:13300} and Figure \ref{fig:18000}) this tendency is accentuated: satellites closer to the warehouse are activated in almost every instance (e.g., 32–33 activations in area A, 26–29 in area B), while the farther from the warehouse the lower the frequency of activation. When switching to TD (distance minimization) the pattern slightly spreads toward central satellites, in more customer-populated areas, confirming that when emissions are not optimized the model tries to reduce both echelons distances reducing the second-echelon intra-route legs. 
Finally, several satellites exhibit very low activation frequencies (often less than 10 times), regardless of the objective or area. These locations are natural candidates for elimination or aggregation, reducing the problem size without materially affecting solution quality. 

\begin{figure}[htbp!]
  \centering
  \resizebox{0.65\linewidth}{!}{%
    \begin{tabular}{@{}cc@{}}
      \begin{subfigure}{0.48\linewidth}
        \centering
        \includegraphics[width=\linewidth]{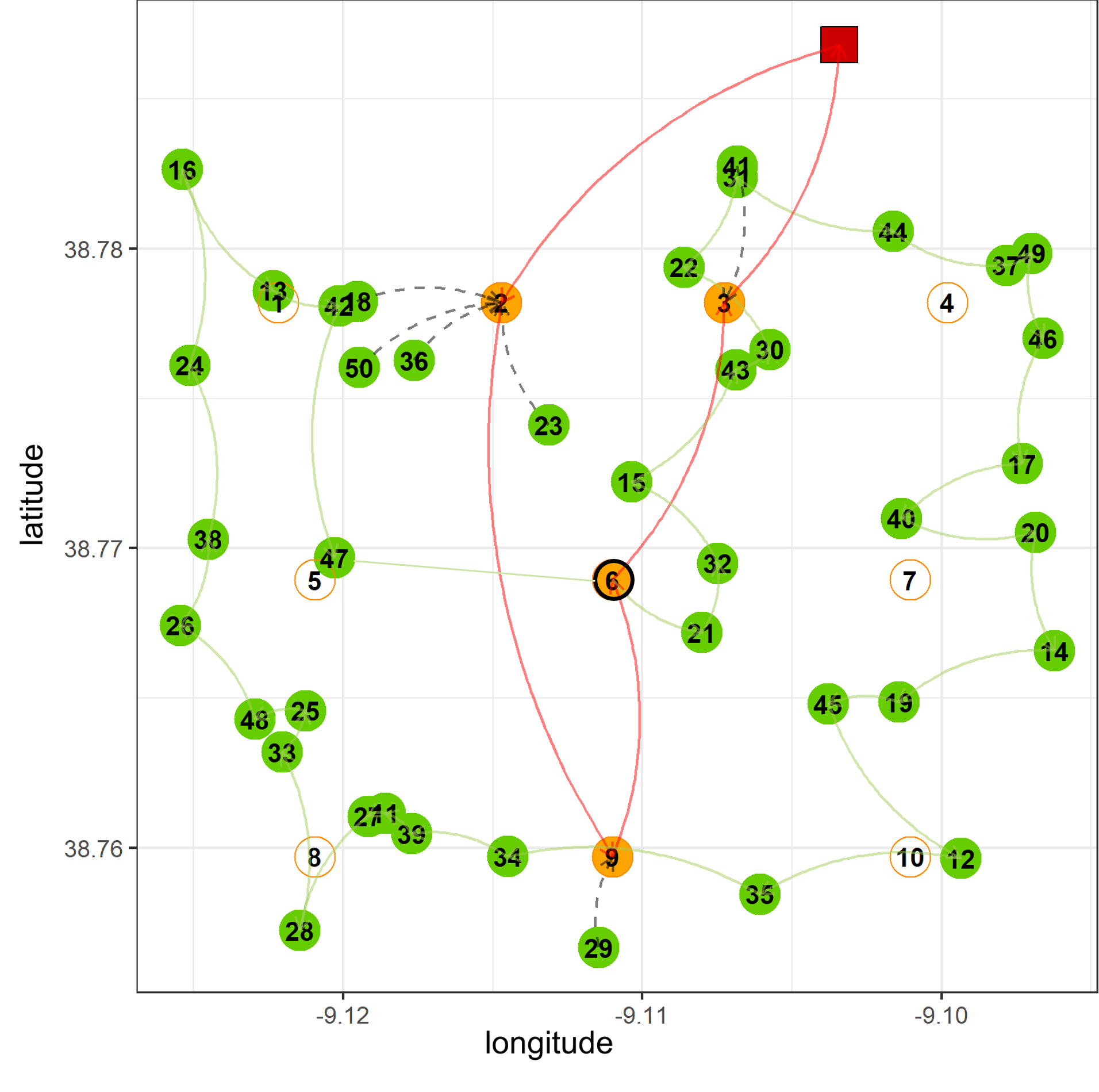}
        \caption{EHC}\label{fig:ex1}
      \end{subfigure} &
      \begin{subfigure}{0.48\linewidth}
        \centering
        \includegraphics[width=\linewidth]{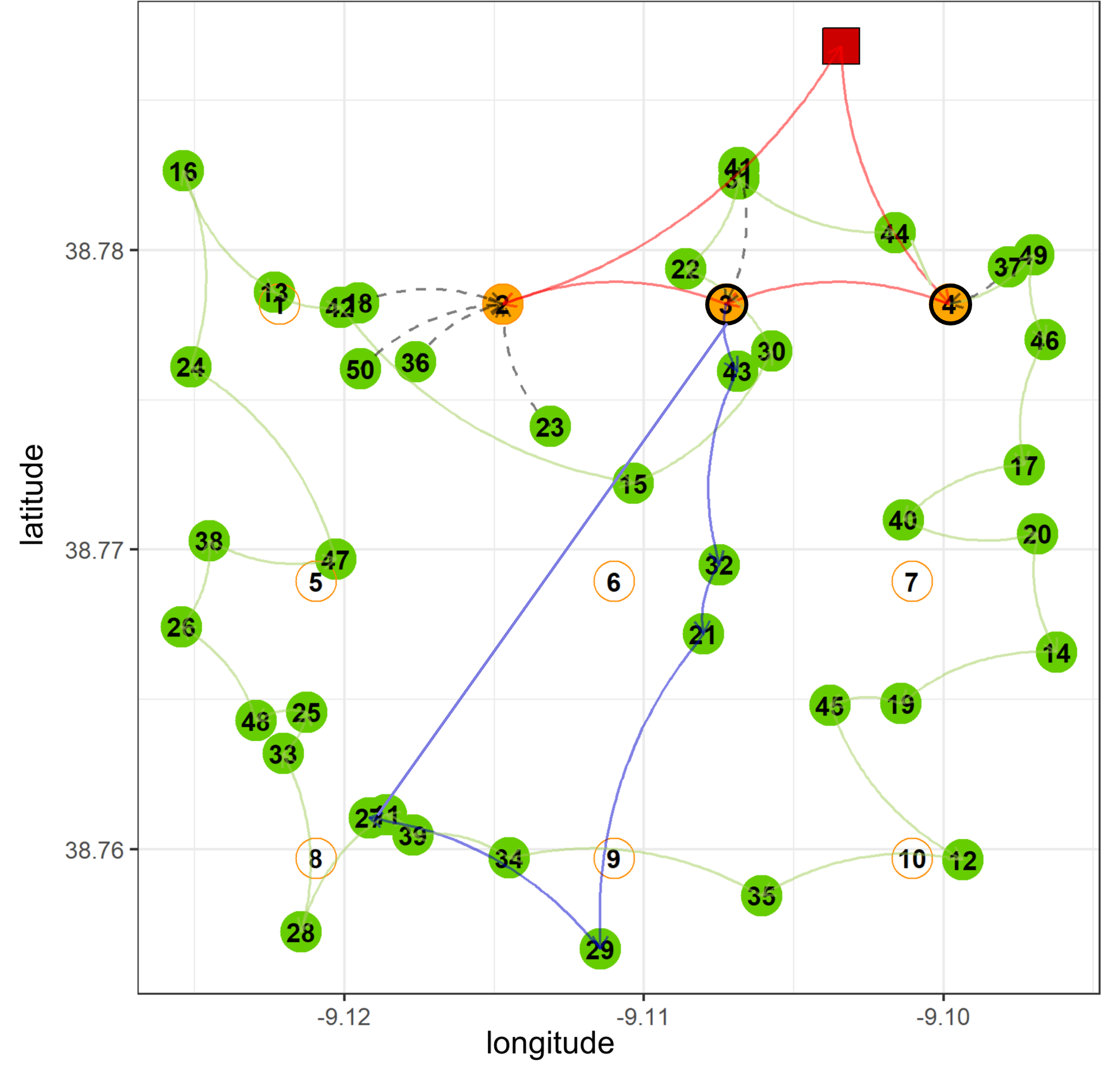}
        \caption{ELC}\label{fig:ex2}
      \end{subfigure} \\
      \begin{subfigure}{0.48\linewidth}
        \centering
        \includegraphics[width=\linewidth]{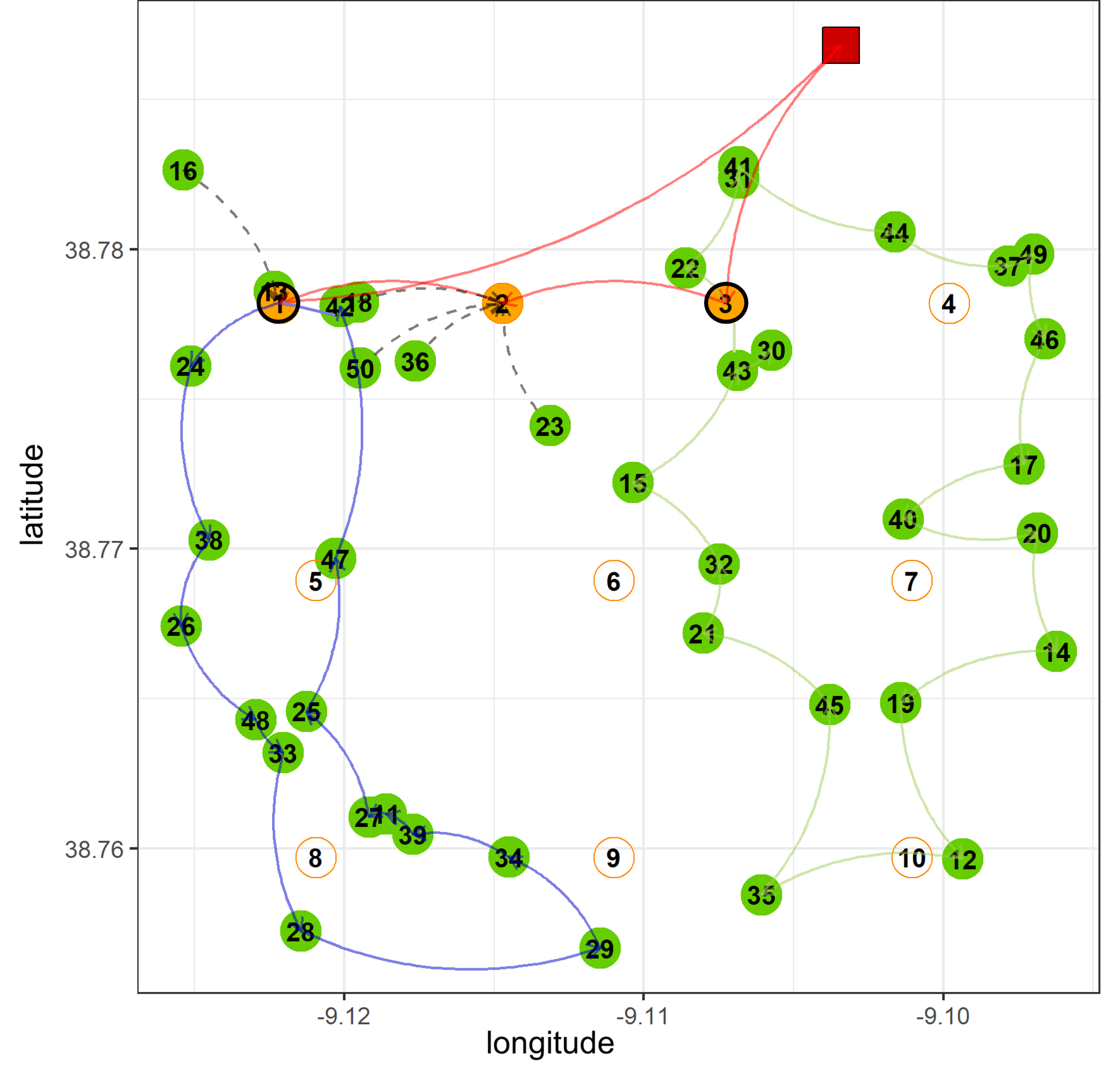}
        \caption{TD}\label{fig:ex3}
      \end{subfigure} &
      \begin{subfigure}{0.48\linewidth}
        \centering
        \includegraphics[width=\linewidth]{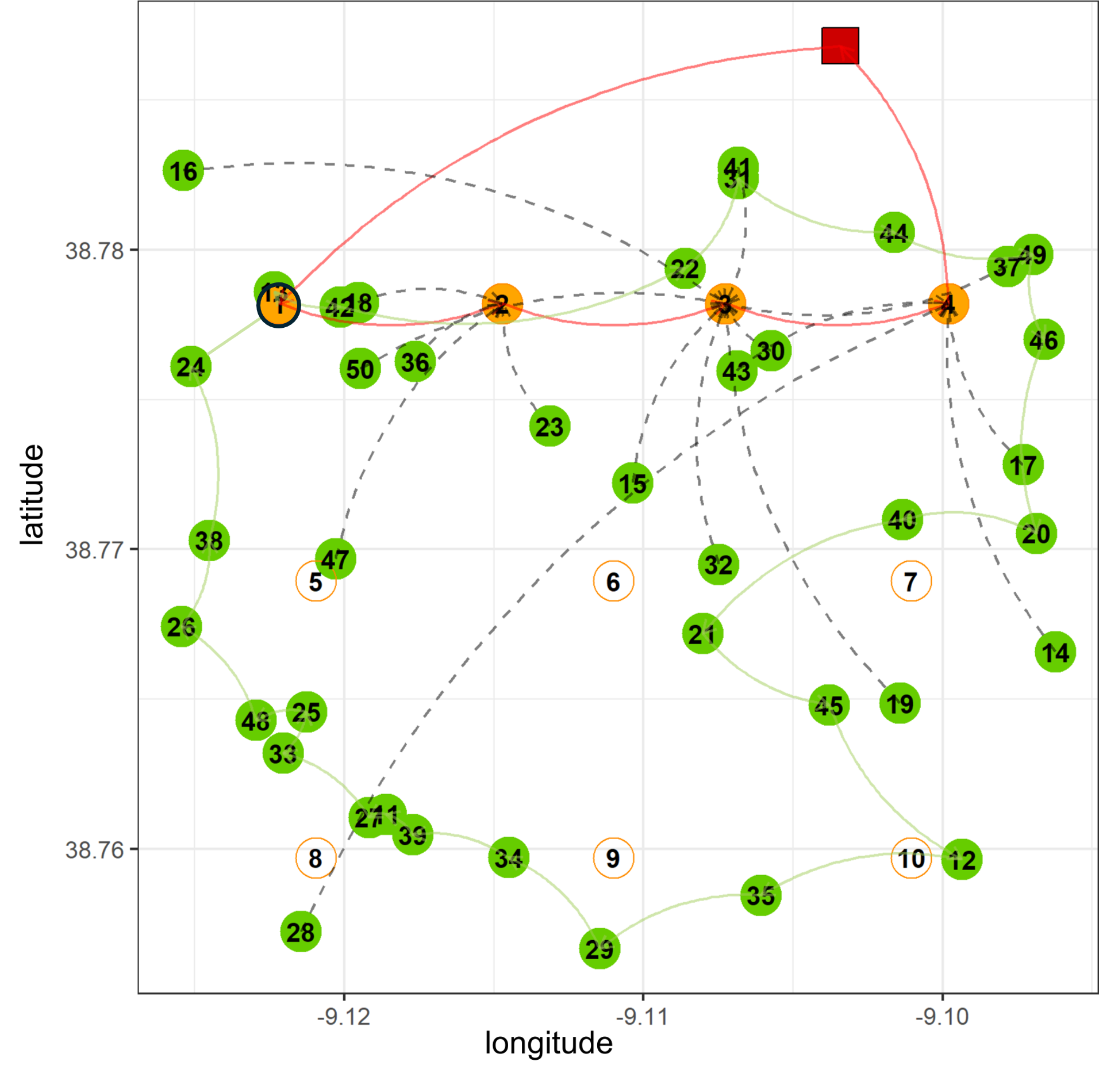}
        \caption{CD}\label{fig:ex4}
      \end{subfigure}
    \end{tabular}%
  }
  \caption{Visual representations of daily activities under different scenarios}
  \label{fig:examples}
\end{figure}

Although satellite–activation numerical statistics do not reveal markedly different patterns across scenarios, the day‑to‑day operations can change substantially. In particular, customers locations, customers behavior, the number and assignment of vehicles, as well as the role of each satellite (activated solely as a pick‑up point versus a home-delivery satellite), may differ. Figure \ref{fig:examples} illustrates these operational contrasts for an instance example with $|C|=30$ customers, $|H|=10$ satellites in area B, highlighting differences that simple aggregate counts cannot capture. In the figure, the warehouse is pictured as a red square. Customers and satellites locations are represented by green and orange dots, respectively. Active satellites are filled in orange and, if a vehicle is assigned they are circled in black. First-echelon routes are represented with red arcs while second-echelon are blue for combustion vehicles and green for zero-emission ones. Customers travels to satellites are drawn as dotted lines.

From Figure \ref{fig:ex1}, representing the minimization of emissions with high capacity vehicles, it can be noted how only zero-emissions vehicles are used. Hence, satellites farther from the warehouse are activated to shorten customers' zero-emission trips (e.g. customer 29 traveling to satellite 9). In contrast, when the capacity of zero-emissions vehicles is reduced as in Figure \ref{fig:ex2}, a combustion vehicle is needed, making it more convenient for the overall emissions to activate only satellites close to the warehouse to compensate for the additional emissions in the second-echelon with a reduction in the first. 
When distances are minimized, as in Figure \ref{fig:ex3} and Figure \ref{fig:ex4}, vehicle emissions are not considered when optimizing routes. In Figure \ref{fig:ex3}, distance optimization leads to a more balanced assignment between the two types of second-echelon vehicles even though the green route could cover the entire demand. This clearly pictures the differences between environmental and operational efficiency. In operational terms, the company benefits from assigning more vehicles, including those with emissions, because this reduces the overall travel distances by splitting the demand across shorter routes. At the same time, whether a customer chooses to pick up an order depends only on the direct saving in distance achieved, without considering the environmental impact. However, environmentally it is prioritized the creation of no-emissions routes for both the company and the customers, even if it results in longer routes. Lastly, when customers are not included in the optimization (Figure \ref{fig:ex4}), company distances are further reduced and satellites close to the warehouse are activated at the expense of longer customer trips that can potentially increase traffic congestion.

\subsection{Customers Behavior Statistics}\label{sec:customers}
In this section, we present some statistics about the customers' behavior. Specifically, we are interested in knowing the percentage of customers picking up their order and the average distance that they travel with and without emissions.
In particular, Table \ref{tab:customers} presents a comparative analysis between the two main scenarios (EHC, which minimizes emissions, and TD, which minimizes total distance) by evaluating the number of customers served at home in column \textit{atHome} and the average distance traveled by those who are required to pick up their orders. For the latter, distances are reported both considering trips with emissions ($dist_{em}$) and trips without emissions ($dist_0$).

\begin{table}[htbp!]
\centering
\scriptsize
\caption{Customers served at home and average customers' traveled distances (in meters) for scenarios EHC and TD}
\begin{tabular}{lrrr|rrr}
\toprule
& \multicolumn{3}{c|}{EHC} & \multicolumn{3}{c}{TD} \\
$|C|$ & atHome (\%) & $dist_{\mathrm{em}}$ & $dist_{0}$ & atHome (\%) & $dist_{\mathrm{em}}$ & $dist_{0}$ \\
\midrule
40  & 36.0 & 738.3 & 373.2 & 37.2 & 421.4 & 290.2 \\
50  & 39.8 & 1039.7 & 321.3 & 46.7 & 357.8 & 274.1 \\
75  & 49.9 & 667.3 & 314.6 & 62.1 & 334.3 & 295.8 \\
100 & 59.2 & 444.2 & 316.6 & 72.7 & 282.5 & 288.4 \\
150 & 71.3 & 385.2 & 326.4 & 87.6 & 329.2 & 319.8 \\
\bottomrule
\end{tabular}
\label{tab:customers}
\end{table}

A clear pattern emerges: in scenario EHC, fewer customers are served at home across all instance sizes compared to scenario TD. This reduction reflects the model's strategy to exploit zero-emission customer travel, especially within allowable green distance bounds. The average $dist_{em}$ under EHC is consistently around or just below 1 km, which is lower than the maximum allowed travel distance of 2 km suggesting that when satellites are too far is suitable for the system to serve the customers at home. At the same time, the corresponding $dist_0$ values typically range between 300 and 370 m, confirming that customer travel is not only sustainable in terms of emissions but also spatially efficient.
In contrast, under scenario TD, a significantly higher share of customers is served directly at home. This configuration prioritizes total distance reduction even at the expense of slightly higher emissions. Nonetheless, the travel burden on pick-up customers remains low in both distance metrics, further suggesting that emissions are indirectly mitigated by minimizing travel distances.
Overall, these results highlight a strategic trade-off: the emissions-oriented approach leverages customer travel more extensively but in a controlled, eco-compliant manner, while the distance-oriented approach favors direct deliveries to reduce total routing length at the system level.

\begin{figure}[htbp]
  \centering
  \begin{subfigure}{0.48\textwidth}
    \centering
    \includegraphics[width=\linewidth]{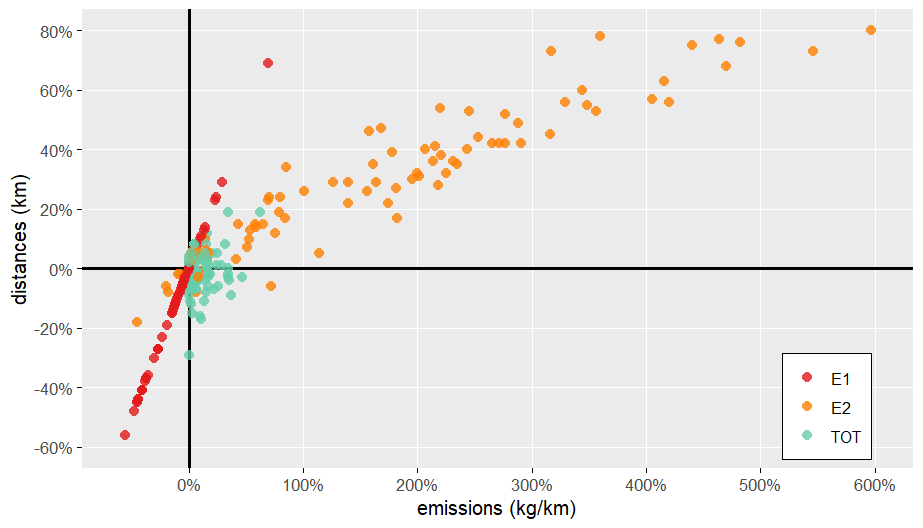} 
    \caption{EHC}
    \label{fig:scatterEHC}
  \end{subfigure}\hfill
  \begin{subfigure}{0.48\textwidth}
    \centering
    \includegraphics[width=\linewidth]{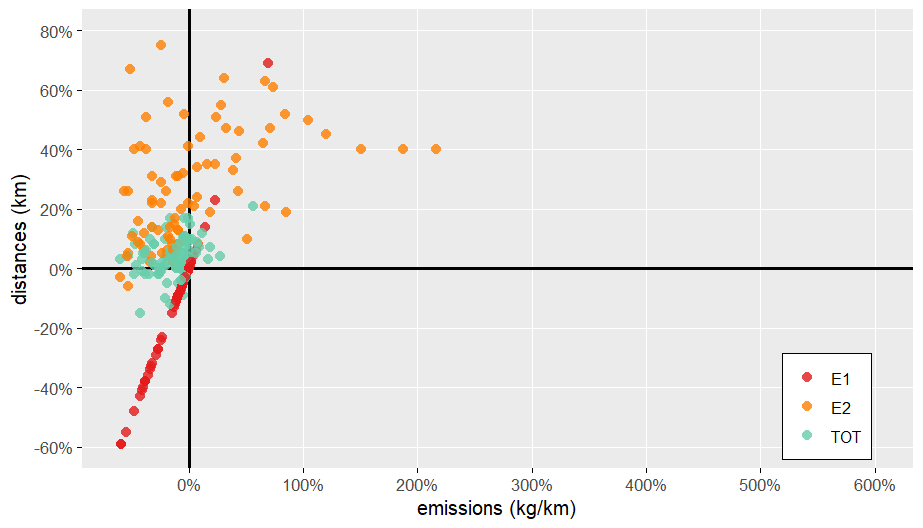} 
    \caption{TD}
    \label{fig:scatterTD}
  \end{subfigure}
  \caption{Percentage variations for EHC and TD}
  \label{fig:scatters}
\end{figure}

Moreover, we re-evaluate all instances by enforcing a full home delivery policy, removing the option for customers to travel to satellites. 
We solve the problem using both the EHC and TD objective functions under this option. The analysis focuses on the resulting variations in total emissions and distances, with a particular emphasis on how the changes manifest across the two-echelons of the delivery network. This allows us to assess the extent to which customer involvement contributes to the efficiency of each optimization strategy and to quantify the additional burden placed on the system when last-mile deliveries are entirely centralized.
Figure \ref{fig:scatterEHC} and Figure \ref{fig:scatterTD} report the percentage variations in emissions (x‑axis) and distances (y‑axis) when the instances are re-solved under a full home‑delivery (HD) system for scenarios EHC and TD, respectively. The variation for any metric is computed as $\frac{f_{HD}-f}{f}*100$ where $f_{HD}$ is the value obtained in the pure home‑delivery setting and $f$ is the corresponding value in the standard model (with customers willing to accept pick-up suggestions).
Results show the variations for the two-echelons separately (dots E1 and E2) and the total variation for the system (TOT).

In Figure \ref{fig:scatterEHC}, the shift to a full home‑delivery system systematically raises total emissions, even when travel distances occasionally decline. This pattern is evident as  the scatter of points lies predominantly in the right half-plane. Although first‑echelon emissions and distances often drop, as shown by the red series E1, they are compensated by the second-echelon growth, with emissions increase reaching up to roughly 600\%, revealing that the burden of home delivery is largely absorbed at this level.

When the objective is total distance minimization (TD), as shown in Figure \ref{fig:scatterTD}, first‑echelon behavior remains similar: emissions and distances in E1 tend to decrease once customer trips are removed. However, the second-echelon and the overall system distances consistently rise, as the company must now cover all last‑mile travels. Moreover, the scatter plot in \ref{fig:scatterTD} shows no clear pattern for emissions variations under TD; the optimization targets operational efficiency regardless of vehicle type, meaning that routes are optimized not considering vehicles emissions. Overall, the figures underscore a common trade-off: if customers are not willing to accept to pick up their orders when needed, this is beneficial to E1, since it allows to activate less satellites but shifts a substantial environmental and operational load onto E2, with different consequences depending on whether emissions or distances is the main objective.

\subsection{Multi-objective Analysis: Pareto Frontiers and Knee-point Identification} \label{sec:multires}

\begin{figure}[htbp]
\centering
\begin{subfigure}[b]{0.48\textwidth}
\centering
\includegraphics[width=\textwidth]{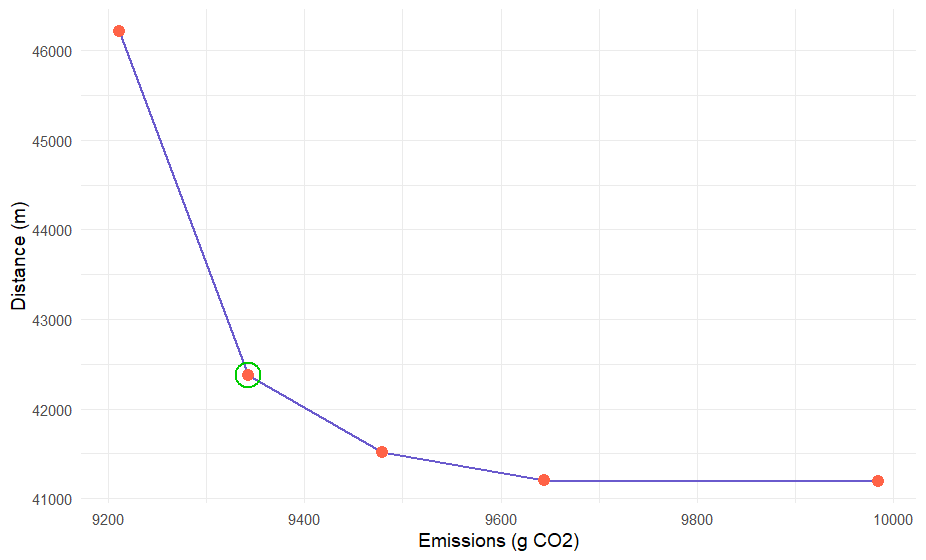}
\caption{50-customer instance}
\label{fig:pareto50}
\end{subfigure}\hfill
\begin{subfigure}[b]{0.48\textwidth}
\centering
\includegraphics[width=\textwidth]{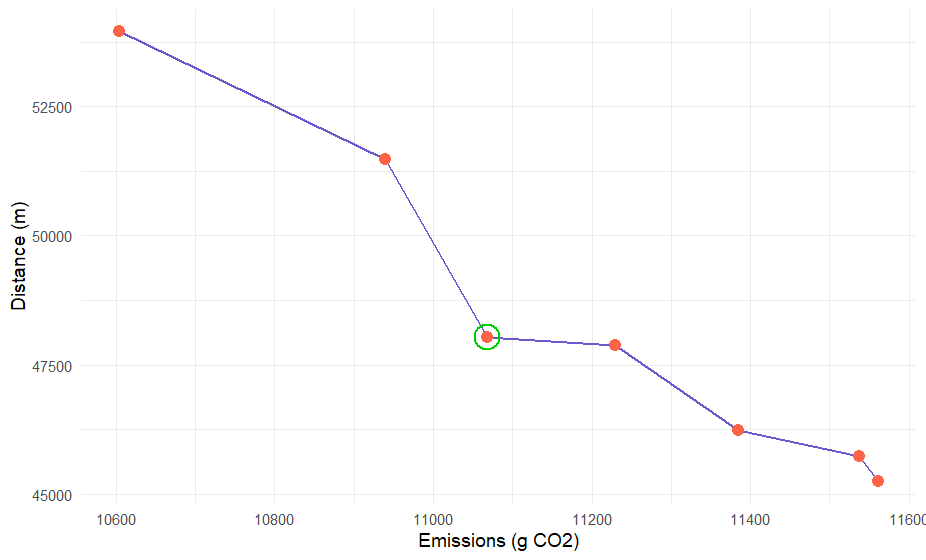}
\caption{75-customer instance}
\label{fig:pareto75}
\end{subfigure}
\caption{Representative approximate Pareto frontiers for two different problem instances.}
\label{fig:pareto}
\end{figure}

To comprehensively evaluate the interplay between environmental sustainability and operational efficiency, we employed a multi-objective optimization approach, specifically the $\epsilon$-constraint method \cite{MesquitaCunha2023}. The primary objective was minimizing emissions, while the distances were progressively constrained within predefined limits, as explained in Section \ref{sec:multi}. 
Figures \ref{fig:pareto50} and \ref{fig:pareto75} illustrate representative approximate Pareto frontiers for 50-customer and 75 customer instances, respectively. Both frontiers exhibit a typical elbow-shaped trend. Initially, a small increase in emissions corresponds to a substantial reduction in traveled distances, highlighting significant operational benefits when allowing slightly higher emissions. This is related to the fact that a combustion vehicle can be inserted in the routes.

The green circle in each figure represents the \textit{knee-point}, computed as detailed in Section \ref{sec:multi}. Identifying this optimal compromise on the Pareto frontier is particularly valuable for practical decision-making. The knee-point corresponds to the solution where the trade-off between the two conflicting objectives, minimization of emissions and distances in our case, becomes most pronounced. At this point, marginal improvements in one objective can only be achieved at a disproportionately high cost in terms of the other, making it the most balanced choice in real-world operations.
For instance, in Figure \ref{fig:pareto50}, the knee-point is located at coordinates $(9343.4,42378.8)$. This solution reveals that attempting to reduce the distance below roughly $42.4$ km leads to steep increases in emissions, meaning operational efficiency gains would come at the expense of environmental sustainability. Conversely, reducing emissions below $9343$ g CO$_2$ would require distance sacrifices that are operationally costly. Thus, this knee-point captures the most stable equilibrium between environmental and economic objectives, making it a highly recommendable solution for practitioners prioritizing both efficiency and sustainability.
Similarly, in Figure \ref{fig:pareto75}, the knee-point is identified at $(11068.3,48041.5)$. Here, the curve is flatter before the knee, indicating that multiple solutions allow simultaneous reductions in both distance and emissions. However, once past the knee-point, any additional improvement in distance efficiency translates into a significant increase in emissions.
From a managerial standpoint, the identification of these knee-points provides a structured way to balance competing goals. Rather than arbitrarily choosing a solution from the Pareto frontier, decision-makers are guided towards the region where the trade-offs are most rational and justifiable. This is particularly important in contexts such as last-mile delivery, where sustainability targets must be reconciled with service efficiency.

\section{Managerial Insights} \label{sec:manager}
The findings of this study yield several practical implications for decision-makers in LMD systems, supporting the definition of sustainable systems. Namely:

\begin{itemize}
    \item[a.] \textbf{Objective choice strongly influences environmental outcomes. (Section \ref{sec:envres})} Optimizing purely for distance minimization can yield modest travel reductions (often below 10\%) but provoke substantial emissions increases (up to 20--30\% in some cases) as it can be seen in Figure \ref{fig:boxplots}. Emission-oriented strategies, conversely, achieve significant sustainability gains with only marginal increases in total distance traveled, indicating that environmental performance should be a primary optimization objective when designing sustainable LMD networks.

    \item[b.] \textbf{Fleet composition is a key lever for sustainability. (Section \ref{sec:envres})} Deploying higher-capacity zero-emission vehicles in the second-echelon reduces the usage on combustion fleets and limits customer travels.
    When low-capacity green vehicles are used, customer pick-up rates increase, shifting part of the delivery burden to customers. This trade-off highlights the importance of aligning vehicle decisions with environmental objectives and customer engagement strategies.

    \item[c.] \textbf{Satellite selection impacts both echelons. (Section \ref{sec:sat})} The model’s results show that satellites closer to the main depot are consistently prioritized to reduce first-echelon emissions . However, when optimizing for distances, more central satellites are activated to shorten second-echelon legs. Operational policies should thus consider dynamic satellite activation rules that balance emissions and service distance targets depending on daily demand patterns.
    \item[d.] \textbf{Incorporating customer behavior can substantially enhance sustainability. (Section \ref{sec:customers})} Encouraging customers to pick up their orders from strategically selected satellites, especially within feasible zero-emission travel distances, leads to notable reductions in second-echelon emissions without significantly increasing the overall traveled distance. Companies can leverage this by promoting voluntary pick-up through incentives.

	\item[e.] \textbf{Multi-objective analysis supports balanced decision-making. (Section \ref{sec:multires})} The $\epsilon$-constraint–based Pareto frontiers provide a quantitative tool to explore the trade-off between environmental and operational objectives. In particular, the identification of \textit{knee-points} highlights the regions where improvements in one objective begin to require disproportionate sacrifices in the other. These points represent practically meaningful configurations where the balance between service efficiency and sustainability is most favorable, guiding decision-makers towards robust and well-justified strategies.
\end{itemize}

Overall, by explicitly considering customer behavior and environmental costs in optimization, companies can achieve a balanced improvement in both environmental sustainability and operational efficiency.

\section{Conclusions} \label{sec: conclusions}
This paper introduces a novel capacitated two-echelon location-routing problem (2E-LRP-ECB) model that explicitly incorporates eco-conscious customer behaviors and heterogeneous vehicles across both echelons, emphasizing environmental sustainability in last-mile delivery operations.

Through extensive computational experimentation, we analyzed various operational strategies, including emissions minimization and distance minimization scenarios, under diverse fleet compositions.
Our findings underscored the importance of integrating customer behavior into last-mile deliveries network design. Requesting customers to pick up orders at satellites significantly enhanced overall sustainability by reducing second-echelon emissions and operational complexity. Moreover, the analyses revealed crucial trade-offs: strategies solely focusing on minimizing traveled distances often resulted in disproportionate increases in emissions, highlighting the necessity of balanced optimization frameworks. Additionally, a multi-objective analysis was proposed providing insightful Pareto frontiers and identified knee-points, effectively guiding decision-makers toward optimal operational-environmental trade-offs. Specifically, the knee-point analysis demonstrated practical benchmarks, illustrating solutions where incremental gains in operational efficiency no longer justify the environmental costs.

Future research could explore incorporating dynamic customer preferences and the integration of additional sustainable transport modes to further enhance model realism and practical applicability. Additionally, extending the study to larger real-world datasets and incorporating stochastic demand scenarios would provide deeper insights into the robustness and adaptability of the proposed framework. 

\section*{Conflict of interest:}
None.
\section*{Acknowledgments}
\begin{itemize}
    \item This work was supported by the European Union’s Horizon Europe programme, Marie Skłodowska-Curie Actions (ERA Fellowship, grant agreement No. 101244242). Views and opinions expressed are however those of the author(s) only and do not necessarily reflect those of the European Union. Neither the European Union nor the granting authority can be held responsible for them.
    \item This work is financed by Portuguese funds through the FCT - Foundation for Science and Technology, I.P., under the projects 2023.13482.PEX and UIDB/00097/2025. 

\end{itemize}

\bibliographystyle{plain} 
\bibliography{main}

\end{document}